\documentclass[11pt]{amsart}
\usepackage[utf8]{inputenc}
\usepackage[a4paper,margin=3cm]{geometry}
\usepackage[pdftex,pdfborder={0 0 0},colorlinks,pagebackref]{hyperref}
\usepackage{latexsym,lmodern,graphicx,enumitem,amssymb,amsmath,amsthm,dsfont}

\newtheorem{thm}{Theorem}[section]
\newtheorem{lem}{Lemma}[section]
\newtheorem{prop}{Proposition}[section]
\newtheorem{cor}{Corollary}[section]

\newtheorem{rem}{Remark}[section]

\numberwithin{equation}{section}

\makeatletter
\def\@MRExtract#1 #2!{#1}
\renewcommand{\MR}[1]{
  \xdef\@MRSTRIP{\@MRExtract#1 !}%
  \href{http://www.ams.org/mathscinet-getitem?mr=\@MRSTRIP}{MR-\@MRSTRIP}}
\makeatother

\newcommand{\R}{\mathbb{R}}
\newcommand{\N}{\mathbb{N}}
\newcommand{\X}{\mathcal{X}}
\newcommand{\T}{\mathbf{T}}
\newcommand{\LSI}{\mathbf{LSI}}

\newcommand{\const}{a}
\newcommand{\cost}{c}

\newcommand{\er}{\mathbb R}

\newcommand{\Cc}{\mathcal C}

\newcommand{\Pp}{\mathcal P}
\newcommand{\Tt}{\mathcal T}

\begin{document}
\title[A variational approach to some transport inequalities]
{A variational approach to some transport inequalities}
\author{Joaquin Fontbona, Nathael Gozlan, Jean-Fran\c cois Jabir}

\thanks{The  first author acknowledges  support  of Basal-Conicyt Center for Mathematical Modeling. The second and third authors were partially supported by the Fondecyt Iniciaci\'on project N 11130705. The first and third authors also thank partial support of  Iniciativa Cientifica Milenio grant  NC130062, Chile. }
\date{March 2016}

\address{Universit\'e Paris Est Marne la Vall\'ee - Laboratoire d'Analyse et de Math\'e\-matiques Appliqu\'ees (UMR CNRS 8050), 5 bd Descartes, 77454 Marne la Vall\'ee Cedex 2, France}
\email{natael.gozlan@u-pem.fr}

\address{Center for Mathematical Modeling,   Universidad de Chile, UMI(2807) UCHILE-CNRS, Casilla 170-3, Correo 3, Santiago, Chile.}
\email{fontbona@dim.uchile.cl}

\address{CIMFAV, Facultad de Ingenieria, Universidad de Valpara\'iso, General Cruz 222, $5^{th}$ floor, Valpara\'iso, Chile.}
\email{jean-francois.jabir@uv.cl}


\maketitle

\begin{abstract}
We relate transport-entropy inequalities to the study of critical points of functionals defined on the
space of probability measures. This approach leads in particular to a new proof of a result by Otto and Villani \cite{OV00} showing that the logarithmic Sobolev inequality implies Talagrand's transport inequality.
\end{abstract}



\section{Introduction}
The aim of this paper is to develop a new variational method for the study of transport-entropy inequalities. This class of inequalities has been introduced by Marton \cite{Mar86,Mar96a, Mar96b}  and Talagrand \cite{Tal96a} in their studies of concentration phenomena for product probability measures. We refer the interested reader to \cite{Vil09, SurveyGL, Led01} for a general exposition on these topics.

The most important transport inequality is certainly the inequality first introduced in \cite{Tal96a} by Talagrand and classically referred to as ``Talagrand's inequality'' or as $\T_2$ inequality in the specialized literature. The inequality $\T_2$ compares two very classical functionals on the space $\mathcal{P}(\X)$ of all probability measures on a given Polish space $\X$: the quadratic Kantorovich distance $W_2(\,\cdot\,,\mu)$ (often called Wasserstein distance) and the relative entropy $H(\,\cdot\,|\mu)$, these two quantities being understood with respect to some fixed reference probability measure $\mu$ on $\X.$ Let us recall the definition of these objects: for all $\nu \in \mathcal{P}(\X)$,
\[
W_2^2(\nu,\mu) = \inf_{\pi} \iint d^2(x,y)\,\pi(dxdy),
\]
where the infimum runs over the set of all couplings $\pi$ between $\nu$ and $\mu$, and
\begin{equation}\label{eq:H}
H(\nu|\mu) = \int \log \left(\frac{d\nu}{d\mu}\right)\,d\nu,
\end{equation}
when $\nu$ is absolutely continuous with respect to $\mu$ (otherwise, one sets $H(\nu|\mu)=+\infty$).\linebreak
A probability measure $\mu$ is said to satisfy the inequality $\T_2(C)$ for some positive constant $C$ if for any probability measure $\nu$ on $\X$,
\[
W_2(\nu,\mu) \leq \sqrt{CH(\nu|\mu)}.
\]
As shown by Talagrand \cite{Tal96a}, the standard Gaussian probability measure on $\R^d$, $d \geq 1$, equipped with the standard Euclidean norm, satisfies $\T_2(2)$.

The inequality $\T_2$, which already has a meaning in terms of comparison of different modes of convergence on the space $\mathcal{P}(\X)$, is also intimately related to the Gaussian concentration of measure phenomenon. Namely, as proved by Talagrand (following a general argument due to Marton) if a probability measure $\mu$ satisfies $\T_2(C)$, then for any positive integer $n$ and for any $1$-Lipschitz function $f$ on $\X^n$ (equipped with the $\ell_2$ distance) it holds
\[
\mu^n(f > m+t) \leq e^{-(t-t_o)^2/C}, \qquad \forall t\geq t_o=\sqrt{C\log(2)},
\]
where $m$ is a median of $f$. This uniform Gaussian control of the tails distributions of Lipschitz maps over product spaces (with constants \emph{independent} on the dimension) - the so-called dimension free concentration property - found numerous applications in various domains (see \cite{Led01, BLM13} for a panorama). This link to concentration of measure is strengthened by the fact that conversely if a probability measure $\mu$ satisfies the property above for some constants $C$ and $t_o$, then it satisfies $\T_2(C)$ (see  \cite{Goz09}).

A natural question is to relate the inequality $\T_2$ to other classical functional inequalities.
An important breakthrough was accomplished in this direction by Otto and Villani \cite{OV00} who first established a clear hierarchy between Talagrand's inequality and the celebrated logarithmic Sobolev inequality. Let us recall the general definition of this well known inequality introduced by Gross \cite{G75}: a probability measure $\mu$ on a metric space $\X$ satisfies the logarithmic Sobolev inequality with a positive constant $C$ - $\LSI(C)$ for short - if for any probability measure $\nu = f\mu$, it holds
\[
H(\nu|\mu) \leq C \int \frac{|\nabla^+ f|^2}{f}\,d\mu,
\]
where in this general context for any function $g:\X \to\R$, and any $x \in \X$, the so-called local slope of $g$ at $x$ is defined by
\begin{equation}\label{eq:local-slope}
|\nabla^+ g|(x) = \limsup_{y \to x} \frac{[
g(y)-g(x)]_+}{d(y,x)}
\end{equation}
(when $x$ is an isolated point in $\X$, then one sets $|\nabla^+ g|(x)=0$).
In \cite{OV00}, Otto and Villani established that, when $\X$ is a smooth connected and complete  Riemannian manifold equipped with its geodesic distance, the logarithmic Sobolev inequality is always stronger than Talagrand's inequality. More precisely, the following holds
\begin{equation}\label{eq:OV}
\LSI(C) \Rightarrow \T_2(4C).
\end{equation}
Roughly speaking, Otto and Villani's proof consists in interpolating $\nu$ and $\mu$ using a certain Fokker-Planck equation (having $\mu$ as limit distribution) and comparing the derivatives of $H$ and $W_2$ along this interpolation. Soon after them, Bobkov, Gentil and Ledoux \cite{BGL01} proposed another proof of the implication \eqref{eq:OV} based on a dual functional formulation of the transport inequality (obtained by Bobkov and G\"otze in \cite{BG99}) and another interpolation technique along this time the solutions of an Hamilton-Jacobi equation. A third proof, based on the characterization of $\T_2$ in terms of concentration discussed above and the well known observation going back to Herbst that $\LSI$ implies Gaussian dimension free concentration (see \textit{e.g.} \cite{Led01}), was proposed by the second author in \cite{Goz09}. It had the advantage over the previous approaches of being immediately generalizable to an abstract metric space framework. It finally turned out that the two other proofs could also be extended to a general metric space context. Namely, Gigli and Ledoux \cite{GL13} have recently adapted the original proof by Otto and Villani to general metric spaces using the general theory of gradient flows as developed in particular in \cite{AGS14}. The proof based on the Hamilton-Jacobi equations has also been adapted to the metric space framework in \cite{GRS14} (improving upon \cite{LV07, Balogh} by removing some unneeded assumptions on the metric measured space appearing in these papers). These three techniques of proof were then re-employed in different settings and for different purposes \cite{Wang04, Wang08, CG06, CGW10, Goz10, GRS11b, GRS13}.

Besides the case of dimension one, where a complete characterization of Talagrand's type inequalities is known (see \cite{Goz07,Goz12}, improving upon \cite{CG06}), the problem of finding sufficient conditions to ensure that a given probability satisfies $\T_2$ is still of great interest (see \cite{Goz10,CGW10} for explicit sufficient conditions on $\R^d$). In this paper, we introduce a new simple method to study transport inequalities and we illustrate it by giving yet another proof of Otto-Villani theorem.

The general idea we develop in the paper is to reduce Talagrand's inequality $\T_2$ (note that the method actually applies to more general transport type inequalities) to the problem of minimizing the function
\begin{equation}
\label{eq:EntW2}
F_a(\nu) = \sqrt{a H(\nu|\mu)} - W_2(\nu,\mu),
\end{equation}
defined for all $\nu \in \mathcal{P}_\mu(\X) := \{\nu \in \mathcal{P}(\X) ; H(\nu|\mu) <\infty\}.$

With this notation in hand we have the following  result.

\begin{lem}\label{lem:argmin} Let $\mu$ be a probability measure on $\X$ and for all $a>0$, denote by $\mathrm{Argmin}(F_a)$ the (possibly empty) set of points $\underline{\nu}$ such that
$F_a(\underline{\nu}) = \inf_{\nu \in \mathcal{P}_\mu(\X)}F_a(\nu).$
\begin{enumerate}
\item The function $F_a$ is bounded from below as soon as $\iint e^{d^2(x,y)/a}\,\mu(dx)\mu(dy)<\infty.$
 \item The probability measure $\mu$ satisfies $\T_2(\const)$ if and only if $\mu\in \mathrm{Argmin}\,(F_\const)$.
 \item The probability measure $\mu$ satisfies $\T_2(\const)$ if and only if for all $a'>a$, $\mathrm{Argmin}\,(F_{a'})=\{\mu\}.$
\end{enumerate}
 \end{lem}
The short proof of this result is postponed to the end of the introduction. Note that the integrability condition given in Item $(1)$ above is not optimal. See Section \ref{Sec:LB} for a discussion and an optimal characterization of the range of parameter $a$ for which $F_a$ is lower bounded in terms of (Gaussian) concentration of measure property.

The question is now to show existence and to characterize minimizers of the function $F_a$.
The existence part is delicate in general, but in the special case where the metric space $(\X,d)$ has a finite diameter, elementary semi-continuity/compactness arguments yield to the conclusion that $\mathrm{Argmin}\,(F_{a})\neq \emptyset$ (see Proposition \ref{prop:attainment} for the finite diameter case and Proposition \ref{prop:attainment2} and Theorem \ref{thm:attainment3} for more general cases). In this introduction, we will always assume that $F_a$ reaches its minimum at (at least) some point, referring to Sections \ref{Sec:Minimization} and \ref{Sec:attainment} for conditions ensuring this property and a thorough discussion of this matter.

In order to state a useful necessary condition satisfied by minimizers of $F_a$, we need to introduce the notion of Kantorovich potentials. According to Kantorovich duality theorem, for all $\nu$,
\[
W_2^2(\nu,\mu) = \sup\left\{ \int \psi\,d\nu + \int \varphi\,d\mu\right\},
\]
where the supremum runs over the set of functions $\psi \in L^1(\nu)$, $\varphi\in L^1(\mu)$ such that $\psi(x) + \varphi(y) \leq d^2(x,y)$ for all $(x,y) \in \X^2$ (see \textit{e.g.} \cite{Vil09}). Under some mild conditions (for instance finite second moments) the supremum is realized by some functions $\overline{\psi},\overline{\varphi}$ related by the following conjugation relations:
\[
\overline{\psi}(x) = \inf_{y}\{d^2(x,y) - \overline{\varphi}(y)\}\qquad \text{and}\qquad \overline{\varphi}(y) = \inf_{x}\{d^2(x,y) - \overline{\psi}(x)\}.
\]
The function $\overline{\psi}$ is usually called a Kantorovich potential for the transport of $\nu$ on $\mu$.

Assuming existence of a minimizer of the function $F_a$ and considering small variations around it, one can prove that a necessary condition for a probability $\underline{\nu} \neq \mu$ to be a minimizer of $F_a$ is to satisfy the following equation
\begin{equation}\label{eq:crit-point}
\lambda \log \left(\frac{d\underline{\nu}}{d\mu} \right) = \overline{\psi} + C,
\end{equation}
where $\lambda = \frac{\sqrt{a}W_2(\underline{\nu},\mu)}{\sqrt{H(\underline{\nu}|\mu)}}$, $C$ is some renormalizing constant and $\overline{\psi}$ is a Kantorovich potential for the transport of $\underline{\nu}$ on $\mu$ (see Theorem \ref{thm:CharactMinimizer} for a general statement).

According to Item (3) of Lemma \ref{lem:argmin}, when a probability measure $\mu$ satisfies $\T_2(a)$ then $\mu$ appears to be the unique minimizer of the functions $F_{a'}$ for $a'>a.$
Therefore a natural sufficient condition to ensure that $\mu$ satisfies $\T_2(a)$ for some $a$ is  to prove that Equation \eqref{eq:crit-point} does not have solution. Indeed, if it is the case, then the only possible minimizer of the function $F_a$ is $\mu$ and so according to Lemma \ref{lem:argmin}, $\mu$ satisfies $\T_2(a).$

Studying the non-linear Equation \eqref{eq:crit-point} appears as a very delicate task. Remarkably, showing that this equation does not have solutions (different to $\mu$) can be achieved easily using the logarithmic-Sobolev inequality. Let us sketch the proof when $\X=\R^d$ is equipped with its usual Euclidean norm. Suppose that $\underline{\nu}\neq \mu$ is a solution of Equation \eqref{eq:crit-point} and that $\mu$ is absolutely continuous with respect to Lebesgue measure. According to a celebrated result by Brenier (see \textit{e.g.} Villani \cite{Vil02}), there exists a transport map $T$ sending $\underline{\nu}$ to $\mu $ (\textit{i.e.} the push forward of $\underline{\nu}$ under the map $T$ is $\mu$) such that
\[
\int |x-T(x)|^2\,\underline{\nu}(dx) = W_2^2(\underline{\nu},\mu).
\]
Moreover, according to classical arguments in optimal transport theory, this map $T$ is related to $\overline{\psi}$ as follows:
\begin{equation}\label{eq:app-transport}
x-T(x) = \frac{1}{2}\nabla \overline{\psi}(x),
\end{equation}
for Lebesgue almost every $x$.
Therefore, reasoning at a formal level, differentiating Equation \eqref{eq:crit-point}, squaring it and integrating with respect to $\underline{\nu}$, yields to
\[
\lambda^2 \int \frac{|\nabla \underline{f}|^2}{\underline{f}^2}\,d\underline{\nu} = 4W_2^2(\underline{\nu},\mu),
\]
denoting by $\underline{f}$ the density of $\underline{\nu}$. By definition of $\lambda$, the latter identity amounts to
\[
a \int \frac{|\nabla \underline{f}|^2}{\underline{f}}\,d\mu = 4H(\underline{\nu}|\mu).
\]
Now, if $\mu$ satisfies $\LSI(C)$ and $a>4C$, this is not possible and so Equation \eqref{eq:crit-point} does not admit solutions, proving that $\mu$ satisfies $\T_2(a)$. This argument is made rigorous in Section \ref{Sec:Applications}, in the general framework of metric spaces, thus giving a new general proof of the implication $\LSI \Rightarrow \T_2$.

In case Equation \eqref{eq:crit-point} admits a non trivial solution $d\underline{\nu} = e^{-V}\,d\mu$, then rewriting \eqref{eq:app-transport} taking into account \eqref{eq:crit-point} one can conclude (at least at a formal level) that
\begin{equation}\label{eq:twistedMM}
d\mu = (\mathrm{Id} + (\lambda/2) \nabla V)_\# \left(e^{-V}\,d\mu\right),
\end{equation}
using the classical notation $S_\#\nu$ to denote the push forward of a measure $\nu$ under a map $S$.
This equation is strongly reminiscent of recent works about the so called \emph{moment measures} \cite{WZ04,D08,BB13,CEK15,L15}. According to a recent result by Cordero-Erausquin and Klartag (see \cite[Theorem 2]{CEK15}), for any probability measure $\mu$ on $\R^d$ having its barycenter at $0$ which is not supported by a lower dimensional subspace of $\R^d$, there exists a convex function $V : \R^d \to \R\cup\{+\infty\}$ essentially continuous (see \cite{CEK15} for a definition) such that $\int e^{-V}\,dx=1$ and
\[
d\mu = \nabla V_\# (e^{-V}\,dx).
\]
Moreover $V$ is unique up to translation. We refer to \cite{K13,K14} for applications of this notion to the study of logarithmically concave probability measures. Our equation \eqref{eq:twistedMM} thus appears as a twisted version of the moment-measure equation above. Actually, the link between our paper and the topic of moment measures is more than formal, since during the preparation of this work, we learned that very recently Santambrogio \cite{S15} recovered the result of \cite{CEK15} through the minimization over the space $\mathcal{P}(\X)$ of a functional $F$ very similar to ours. Some further details are provided in Section \ref{Sec:Moment-Measures}.

The rest of the paper is organized as follows.
Section \ref{Sec:Minimization} is devoted to the study of a class of functionals generalizing \eqref{eq:EntW2}. These functionals are of the form $F_a(\nu) = \alpha(aH(\nu|\mu)) - \beta\left(\mathcal{T}_c(\nu,\mu)\right)$, $\nu \in \mathcal{P}_\mu(\X)$, where $\alpha$ and $\beta$ are given functions on $[0,\infty)$ and $\mathcal{T}_c$ is an optimal transport cost associated to some general cost function $c$ on $\X$. The question of existence of a minimizer is discussed and an equation generalizing \eqref{eq:crit-point} is derived for those minimizers.
Section \ref{Sec:Applications} is dedicated to applications. We prove in particular different variants of the Otto-Villani Theorem in metric spaces.
In Section \ref{Sec:attainment}, we prove Theorem \ref{thm:attainment3} (stated in Section \ref{Sec:Minimization}) establishing the existence of a minimizer for functionals of the form $F_a(\nu) = a H(\nu|\mu) - \mathcal{T}_c(\nu,\mu)$, $\nu \in \mathcal{P}_\mu(\X)$, under a weak (and actually minimal) concentration of measure assumption for $\mu$.
Finally Section \ref{Sec:Moment-Measures} contains remarks about the links between our work and the main results of \cite{CEK15} and \cite{S15}.

\proof[Proof of Lemma \ref{lem:argmin}]
(1) This is a consequence of the general result proved in Proposition \ref{prop:bounded-below}.
(2) The probability $\mu$ satisfies $\T_2(\const)$ if and only if $F_\const(\nu)\geq 0$ for all $\nu\in \mathcal{P}(\X)$. Since $F_\const(\mu)=0$, this is equivalent to the condition $\mu\in \mathrm{Argmin}\,(F_\const).$ (3) Suppose that $\mu$ satisfies $\T_2(\const)$ for some $\const>0$. Then, if $\const'>\const$ then $\mu$ also satisfies $\T_2(\const')$. If $\underline{\nu}\in \mathrm{Argmin}\,(F_{\const'})$, it holds $\const'H(\underline{\nu}|\mu)=W_2^2(\underline{\nu},\mu)\leq \const H(\underline{\nu}|\mu)$. Since $\const'>\const$, the only possibility is that $\underline{\nu}=\mu$. The converse is immediate.
\endproof

\section{Minimization of a class of functionals on the space $\mathcal{P}(\X)$}\label{Sec:Minimization}
In this section we introduce a class of functionals involving the relative entropy and a optimal general transport cost. Then we study the minimization problem of these functionals and we give a characterization of the optimizers.

We recall that in all the paper $(\X,d)$ is a Polish space and that $\mathcal{P}(\X)$ denotes the set of all Borel probability measures on $\X.$
\subsection{Definitions}\label{Sec:Def}
First let us recall the definition of optimal transport costs. Given a cost function
\[
\cost:\X\times\X \rightarrow \er^+
\]
that we will assume hereafter to be continuous, for all probability measures $\nu_1,\nu_2$ on $\X$, one denotes by $\mathcal{T}_\cost(\nu_1,\nu_2)$ the optimal transport cost between $\nu_1$ and $\nu_2$ defined by
\[
\mathcal{T}_\cost(\nu_1,\nu_2)=\inf \iint \cost(x,y)\,\pi(dxdy)\in [0,+\infty],
\]
where the infimum runs over the set of all couplings $\pi\in \mathcal{P}(\X^2)$ having $\nu_1$ and $\nu_2$ as first and second marginal distributions.

Given a probability measure $\mu$ on a $\X$, one denotes by $\mathcal{P}_\mu(\X)$ the set of probability measures $\nu$ such that $H(\nu|\mu)<\infty$ (recall the definition \eqref{eq:H} of the relative entropy functional). Then we consider the functional
$F_\const:\mathcal{P}_\mu(\X)\to \R\cup\{- \infty\}$ where $\const>0$ defined by
 \begin{equation}\label{eq:Fa}
 F_\const(\nu)=\alpha\left(\const H(\nu|\mu)\right)-\beta\left(\mathcal{T}_\cost(\nu,\mu)\right),\quad \forall \nu \in \mathcal{P}_\mu(\X),
 \end{equation}
where $\alpha,\beta:\er^+ \rightarrow \er^+$ are two $\Cc^1$-functions on $(0,+\infty)$ and $\beta$ is assumed to be non-decreasing. In all what follows we will also assume that the functions $\alpha,\beta$ are such that for all $b\geq 0$ the function
\[
t\mapsto \alpha(t)-\beta(t+b)
\]
is bounded from below on $(0,\infty)$ and for all $\lambda>1$,
\[
\alpha(\lambda t) - \alpha(t) \to +\infty, \text{ when } t \to +\infty.
\]
The functional considered in the introduction corresponds to $\alpha(t) = \beta(t) = \sqrt{t}$ and $c(x,y)=d^2(x,y)$. Other choices will be considered in Section \ref{Sec:Applications}.

Throughout the paper, we will often deal with a particular class of cost functions generalizing power cost functions of the form $c(x,y) = d^p(x,y)$, $x,y\in \X$, $p\geq 1$.
This class of cost functions is introduced in the following lemma.
\begin{lem}\label{lem:power-type}
Let $c:\X \times \X \to \R^+$ be a cost function of the form
\[
c(x,y) = \phi(d(x,y)),\qquad x,y \in \X
\]
where $\phi : \R^+ \to \R^+$ is some convex function such that $\phi(0)=0$, $\phi(x)>0$ if $x>0$ and $\sup_{x >0} \phi(2x)/\phi(x) <\infty.$
Define $p_o = \sup_{x>0} x \phi'(x) / \phi(x)$, where $\phi'$ denotes the right-derivative of $\phi.$ Then $1\leq p_o <\infty$, and the function $\tilde{d}$ defined by
\begin{equation}\label{eq:dtilde}
\tilde{d}(x,y) = c^{1/p_o}(x,y),\qquad x,y\in \X
\end{equation}
is a distance on $\X$ inducing the same topology as $d.$
\end{lem}

In all what follows, a cost function as in Lemma \ref{lem:power-type} will be referred to as a \emph{power type} cost function, and the number $p_o$ associated to it will be called its \emph{exponent}.

\proof
Define $K=\sup _{x >0} \phi(2x)/\phi(x)< \infty$. By convexity, $\phi(x)/x \leq \phi'(x)$ and so $p_0 \geq 1.$ On the other hand, $\phi(2x) \geq \phi(x) + \phi'(x) x$ and so $p_0\leq K-1<\infty.$  Set $\omega (x) =\phi^{1/p_o}(x)$, $x \geq0$. Then
\[
\frac{d}{dx}\left\{\frac{\omega(x)}{x}\right\} =\frac{1}{x^2}\left[\frac{1}{p_o} \phi^{1/p_o}(x) \frac{x\phi'(x)}{\phi(x)} - \phi^{1/p_o}(x)\right] \leq  0.
\]
Therefore the function $x \mapsto \omega(x)/x$ is non-increasing. As a result, the function $\omega$ is sub-additive: $\omega(a+b) \leq \omega(a) + \omega(b)$, for all $a,b \geq0.$
One concludes from this that the function $\tilde{d} = c^{1/p_o}$ is a distance on $\X$. The last assertion follows from the fact that $d$ and $\tilde{d}$ define the same set of converging sequences.
\endproof

\subsection{Conditions for lower boundedness}\label{Sec:LB}
Our general purpose being to study the minimization of $F_a$ over $\mathcal{P}_\mu(\X)$, let us begin with the following simple observation showing that $F_a$ is bounded from below when the cost $c$ is exponentially integrable with respect to $\mu\otimes \mu.$

\begin{prop}\label{prop:bounded-below}
Let $c:\X^2 \to \R^+$ be a continuous cost function and $\mu \in \mathcal{P}(\X)$ be such that
\begin{equation}\label{eq:expmoment}
I_\delta:=\iint e^{\delta c(x,y)}\,\mu(dx)\mu(dy)<\infty,
\end{equation}
for some $\delta>0$.
Then, for all $\nu \in \mathcal{P}(\X)$
\[
\iint c(x,y)\,\nu(dx)\mu(dy) \leq \frac{1}{\delta}H(\nu|\mu) + \frac{e^{-1}}{\delta} I_\delta.
\]
In particular, for all $a \geq 1/\delta$, the function $F_a$ is bounded from below on $\mathcal{P}_\mu(\X)$. Moreover, for all $a >1/\delta$, the level sets $\{F_a \leq r\}$, $r\in \R$, of the function $F_a$ are precompact for the weak topology on $\Pp(\X)$.
\end{prop}
\proof
Consider the functions $\theta(x) = x\log(x)$, $x >0$, and $\theta^*(y) = \sup_{x>0} \{xy- \theta(x)\}$, $y\in \R$. An easy calculation shows that $\theta^*(y)= e^{y-1}$ for all $y\in \R$.
Assuming without loss of generality that $\nu \in \mathcal{P}_\mu(\X)$ and using the immediate Young's inequality
\begin{equation}\label{eq:Young}
xy \leq \frac{1}{\delta}\theta(x) + \frac{1}{\delta}\theta^*(\delta y),\qquad  x >0, y\in \R,
\end{equation}
yields to
\begin{align*}
\iint c(x,y)\,\nu(dx)\mu(dy)&= \iint c(x,y)\frac{d\nu}{d\mu}(x)\,\mu(dx)\mu(dy)\\
& \leq \frac{1}{\delta}\iint \theta\left(\frac{d\nu}{d\mu}(x)\right)\,\mu(dx)\mu(dy)
+\frac{1}{\delta}\iint \theta^*(\delta c(x,y))\,\mu(dx)\mu(dy)\\
& = \frac{1}{\delta} H(\nu|\mu) + \frac{e^{-1}}{\delta} I_{\delta}.
\end{align*}
Since $\Tt_\cost(\nu,\mu) \leq \iint c(x,y)\,\nu(dx)\mu(dy)$ and $\beta$ is non-decreasing one concludes that
\[
F_a(\nu) = \alpha(a H(\nu|\mu)) - \beta(\mathcal{T}_c(\nu,\mu)) \geq \alpha(aH(\nu|\mu)) - \beta(\delta^{-1}H(\nu|\mu) + b),
\]
with $b=e^{-1}I_\delta/\delta.$
Writing, for all $h\geq0$,
\[
\eta(h):=\alpha(ah)-\beta(\delta^{-1}h+b) = \alpha(ah)-\alpha(\delta^{-1}h) + \alpha(\delta^{-1}h) - \beta(\delta^{-1}h+b),
\]
one sees using the assumptions made on $\alpha$ and $\beta$ that $\eta(h)$ is bounded from below.
This implies that $F_a$ is bounded from below on $\mathcal{P}_\mu(\X)$. Moreover, if $a>\delta^{-1}$, then $\eta(h) \to \infty$, as $h \to +\infty.$ Therefore, for each $r\in \R$, there exists $h_r\geq 0$ such that
\begin{equation}\label{eq:level-sets}
\{\nu \in \mathcal{P}_\mu(\X) : F_a(\nu) \leq r\} \subset \{\nu \in \mathcal{P}(\X) :  H(\nu|\mu) \leq h_r\}.
\end{equation}
According to \textit{e.g.} \cite[Lemma 6.2.12]{DZ98}, the level sets of $\nu\rightarrow H(\nu|\mu)$ are compact for the weak topology of $\Pp(\X)$. This completes the proof.
\endproof

\begin{rem}The use of integrability conditions to prove transport inequalities is now very classical. Let us mention in particular the seminal paper by Djellout-Guillin-Wu \cite{DGW04} establishing the equivalence between the transport inequality $\T_1$ and a Gaussian integrability condition. This approach was then further developed by Bolley-Villani \cite{BV05} and the second named author \cite{Goz06}. A slightly different point of view was proposed by Milman in the paper \cite{Mil12}, where a quantitative equivalence is established between the concentration of measure properties of a measure and transport inequalities involving the Kantorovich $W_1$ distance. This point of view was then further developed in a paper by Roberto, Samson and the second named author  \cite{GRS11a} in terms of non-tight transport inequalities involving the $W_2$ distance. This last result is recalled in Proposition \ref{prop:bounded-below-refined} below.
\end{rem}

According to Proposition \ref{prop:bounded-below} above, $F_a$ is bounded from below on $\mathcal{P}_\mu(\X)$ as soon as the integrability condition \eqref{eq:expmoment} is fulfilled for $\delta=1/a$. It turns out that this integrability condition is not necessary to ensure that $F_a$ is bounded from below. Indeed, let $\X = \R^d$ equipped with its usual Euclidean distance and let $\mu$ be the standard Gaussian measure. On the one hand, it is easily seen that $I_\delta = \iint e^{\delta d^2(x,y)}\,\mu(dx)\mu(dy)$ is finite only when $\delta < 1/4.$ Therefore, according to Proposition \ref{prop:bounded-below}, the functional $F_a$ is bounded from below for all $a>4.$ But, on the other hand, since $\mu$ satisfies Talagrand's inequality $\T_2(2)$ (see \cite{Tal96a}), we see that $F_a\geq 0$ on $\mathcal{P}_\mu(\X)$ for all $a \geq 2$. So for $a \in [2,4]$, $F_a$ is bounded from below and $I_{1/a}=+\infty.$

Fortunately, for power type cost functions (as defined in Lemma \ref{lem:power-type}) Proposition \ref{prop:bounded-below} can be improved. Indeed, when $c$ is a power type cost function (and $\alpha = \beta = \mathrm{Id}$) the range of parameters $a$ for which the function $F_a$ is bounded from below can be completely determined in terms of a concentration of measure inequality for $\mu$ (with respect to the metric $\tilde{d}$ introduced in Lemma \ref{lem:power-type}).

\begin{prop}\label{prop:bounded-below-refined}
Let $c:\X \times \X \to \R^+$ be a power type cost function and $p_o$ be its exponent. For all $a>0$, consider the function $F_a$ defined by $F_a(\nu)=aH(\nu|\mu) - \mathcal{T}_c(\nu,\mu)$, $\nu \in \mathcal{P}_\mu(\X)$.

The following propositions are equivalent:
\begin{enumerate}
\item There exist $a>0$, $b\geq0$ such that the function $F_a\geq-b$ on $\mathcal{P}_\mu(\X)$.
\item There exist $a'>0$ and $r_o \geq0$ such that the probability measure $\mu$ satisfies the following concentration of measure property: for all $A \subset \X$ such that $\mu(A) \geq 1/2$, it holds
\begin{equation}\label{eq:concentration}
\mu(A_r) \geq 1-e^{-(r-r_o)^{p_o}/a'},\qquad \forall r\geq r_o,
\end{equation}
where $A_r = \{ x \in \X: \exists y \in A, \tilde{d}(x,y) \leq r \}$ and $\tilde{d}(x,y) = c^{1/p_o}(x,y)$, $x,y \in \X$.
\end{enumerate}
More precisely,
\begin{itemize}
\item (1) $\Rightarrow$ (2) with $a'=a$, and $r_o=(a\log(2))^{1/p_o} + 2 b^{1/p_o}$,

and

\item (2) $\Rightarrow$ (1) for all $a=ta'$ with $t>1$ and for some $b$ depending on $r_o,a'$ and $t$.
\end{itemize}
Moreover, assuming (2), the level sets $\{F_a\leq r\}$, $a >a'$, $r\in \R$ are precompact  for the weak topology on $\mathcal{P}(\X).$
\end{prop}
To summarize, for the functionals $F_a$ considered in Proposition \ref{prop:bounded-below-refined}, it holds
\begin{equation}\label{eq:exact-range}
\inf\{ a >0 : F_a \text{ is bounded below}\} = \inf\{ a '>0 : \eqref{eq:concentration} \text{ holds for some $r_o\geq0$}\}.
\end{equation}

\proof
The implication (1) $\Rightarrow$ (2) follows a well known general argument due to Marton \cite{Mar86} briefly sketched below. Using Lemma \ref{lem:power-type} and $c(x,y) = \tilde{d}^{p_o}(x,y)$, $x,y\in \X$, one can represent $\mathcal{T}_c$ as a Wasserstein distance: for all $\nu_1,\nu_2 \in \mathcal{P}(\X)$,
\[
\mathcal{T}_c(\nu_1,\nu_2) = \widetilde{W}_{p_o}^{p_o}(\nu_1,\nu_2) := \inf_{(X,Y)} \mathbb{E}\left[\tilde{d}^{p_o}(X,Y)\right],
\]
where the infimum runs over the set of couples $(X,Y)$ of random variables such that $X$ has law $\nu_1$ and $Y$ has law $\nu_2.$ Let $A \subset \X$ be a Borel set such that $\mu(A)\geq 1/2$ and define $B=\X \setminus A_r$, for some $r\geq0$, and consider the probability measures $d\mu_A = \frac{\mathbf{1}_A}{\mu(A)} \,d\mu$ and $d\mu_B = \frac{\mathbf{1}_B}{\mu(B)} \,d\mu$. Using the triangle inequality for $W_p$, we get
\begin{align*}
 \widetilde{W}_{p_o}(\mu_A,\mu_B) &\leq  \widetilde{W}_{p_o}(\mu_A,\mu) +  \widetilde{W}_{p_o}(\mu,\mu_B)\\
& \leq \left(aH(\mu_A|\mu)\right)^{1/p_o} + \left(aH(\mu_B|\mu)\right)^{1/p_o} + 2 b^{1/p_o}\\
& \leq \left(a \log(2)\right)^{1/p_o}+b^{1/p_o} + \left(a\log(1/\mu(B)) \right)^{1/p_o},
\end{align*}
where the second line comes from the inequality $ \widetilde{W}_p(\nu,\mu) \leq \left(aH(\nu|\mu)\right)^{1/p_o} + b^{1/p_o}$, which is easily deduced from the assumption that $F_a \geq -b$. Using that $W_p(\mu_A,\mu_B) \geq r$  (which follows at once from the definition of the set $B$) one easily gets \eqref{eq:concentration}.

The implication (2) $\Rightarrow$ (1) is adapted from \cite[Corollary 2.20]{GRS11a}. More precisely, one easily derives from \eqref{eq:concentration} that for all $a'' > a'$, there is some $M_{a''}\geq 1$ such that for all $A \subset \X$ with $\mu(A) \geq 1/2$, it holds
\[
\mu(A_r) \geq 1-M_{a''}e^{-r^{p_o}/a''},\qquad \forall r\geq 0.
\]
If $p_o =2$ then Corollary 2.20 of \cite{GRS11a} shows that for all $a >a''$, there is some $b$ depending on $a'', M_{a''}$ and $t=a/a''$ such that
\[
\widetilde{W}_2^2(\nu,\mu) \leq aH(\nu|\mu) + b,\qquad \forall \nu\in \mathcal{P}_\mu(\X),
\]
which in other words means that $F_a\geq -b$ on $\mathcal{P}_\mu(\X).$ This proves the implication (2) $\Rightarrow$ (1) in the case $p_o=2$. Now if $p_o \neq 2$, it turns out that the proof of Corollary 2.20 can be very easily adapted (just replacing $W_2$ by $W_p$, $r^2$ by $r^p$, etc\ldots) yielding exactly as before to the conclusion that $F_a$ is bounded from below on $\mathcal{P}_\mu(\X)$ for all $a>a'.$

It remains to show that if \eqref{eq:concentration} holds, then the level sets $\{F_a \leq r\}$, $r\in \R$ are precompact for $a>a'.$ Indeed, for all $\nu \in \mathcal{P}_\mu(\X)$, it holds
\[
F_a(\nu) = F_{(a+a')/2}(\nu) + \frac{(a-a')}{2}H(\nu|\mu) \geq m + \frac{(a-a')}{2}H(\nu|\mu),
\]
where $m = \inf_{\nu \in \mathcal{P}_\mu(\X)} F_{(a+a')/2}(\nu)$ which is finite according to the implication (2) $\Rightarrow$ (1). From this follows that $\{F_a \leq r\} \subset \{H(\,\cdot\,|\mu) \leq 2(r-m)/(a-a')\}$. Since the level sets of the relative entropy are compact for the weak topology, the proof is complete.
\endproof

To conclude this section, let us give an elementary comparison between the integrability and concentration conditions considered above.
\begin{lem}\label{lem:comparison}
Let $c:\X \times \X \to \R^+$ be a power type cost function with exponent $p_o$ and associated distance $\tilde{d}=c^{1/p_o}.$
\begin{enumerate}
\item If $I_\delta:=\iint e^{\delta c(x,y)}\,\mu(dx)\mu(dy)<\infty$ for some $\delta>0$, then $\mu$ satisfies the concentration inequality \eqref{eq:concentration} with $a'=1/\delta$ and $r_o = (\delta \log(2I_\delta))^{1/p_o}$.
\item If $\mu$ satisfies the concentration inequality \eqref{eq:concentration} for some $a'>0$ and $r_o\geq 0$, then for all $x_o \in \X$, $\int e^{\delta c(x_o,y)}\,\mu(dy) <\infty$ for all $\delta<1/a'$ and $\iint e^{\delta c(x,y)}\,\mu(dx)\mu(dy) <\infty$ for all $\delta<1/(a'2^{p_o-1}).$
\end{enumerate}
\end{lem}
\proof
(1) Let $A$ be such that $\mu(A) \geq 1/2$ and set $B= \X \setminus A_r$; then it holds
\[
I_\delta=\iint e^{\delta c(x,y)}\,\mu(dx)\mu(dy) \geq \int_A \int_B e^{\delta c(x,y)}\,\mu(dx)\mu(dy)\geq e^{\delta r^{p_o}}\mu(A)\mu(B).
\]
Writing that $\mu(A)\geq 1/2$, the announced inequality easily follows.
(2) Let $m$ denote the median of the function $x\mapsto \tilde{d}(x,x_o)$ where $x_o$ is some arbitrary point in $\X.$ According to the classical formulation of concentration of measure in terms of deviation inequalities for $1$-Lipschitz functions (see \textit{e.g.} \cite{Led01}), it holds
\[
\mu( \{y \in \X :  \tilde{d}(x_o,y)>m+r\}) \leq e^{-[r-r_o]_+^{p_o}/a'},\qquad \forall r\geq 0,
\]
where $[\,\cdot\,]_+$ denotes the positive part function.
By an integration by part
\begin{align*}
\int e^{\delta c(x_o,y)}\,\mu(dy) = \int e^{\delta \tilde{d}^{p_o}(x_o,y)}\,\mu(dy) &= 1+ \delta \int_0^{+\infty} e^{\delta v}\mu(\{y : \tilde{d}(x_o,y)> v^{1/p_o}\})\,dv\\
& \leq 1+ \delta \int_0^{+\infty} e^{\delta v}e^{-[v^{1/p_o}-r_o-m]_+^{p_o}/a'}\,dv
\end{align*}
and this last integral is clearly finite if and only if $\delta <1/a'.$ Now using the inequality
\[
(u+v)^{p_o} \leq 2^{p_o-1}u^{p_o} + 2^{p_o-1}v^{p_o},\qquad u,v \geq0
\]
and the triangle inequality for $\tilde{d}$, one concludes that $\iint e^{\delta c(x,y)}\,\mu(dx)\mu(dy) <\infty$ as soon as $2^{p_o-1}\delta <1/a'$, which completes the proof.
\endproof

\begin{rem}\label{rem:comparison}
For a given probability measure $\mu$ and a given power type cost function $c$, let us denote by
\begin{itemize}
\item $a_o$ the infimum of the $a>0$ such that $F_a=aH(\,\cdot\,|\mu)-\mathcal{T}_c(\,\cdot\,,\mu)$ is bounded below on $\mathcal{P}_\mu(\X)$,
\item $a'_o$ the infimum of the $a'>0$ such that $\mu$ satisfies the concentration inequality \eqref{eq:concentration} for some $r_o\geq0$,
\item $\delta_o$ the supremum of the $\delta >0$ such that $\iint e^{\delta c(x,y)}\,\mu(dx)\mu(dy) <\infty$,
\item $\delta_o'$ the supremum of the $\delta' >0$ such that $\int e^{\delta c(x_o,y)}\,\mu(dy) <\infty$ for all $x_o \in \X.$
\end{itemize}
We established in \eqref{eq:exact-range} that $a_o=a_o'$. We can complete this result by the following inequalities immediately deduced from Lemma \ref{lem:comparison}:
\[
1/(2^{p_o-1}\delta_o) \leq 1/\delta_o' \leq a_o' \leq 1/\delta_o.
\]
The example of the standard Gaussian discussed above showed that $a_o$ can be strictly less than $1/\delta_o$. We do not know if $a_o'=1/\delta_o'$.
\end{rem}

\subsection{Attainment of the minimum}
We begin with a simple result showing that, when the cost function is continuous and bounded, the function $F_a$ attains its minimal value. Note that in the case of a bounded cost function $c$, the definition of the function $F_a$ makes sense over the whole $\mathcal{P}(\X).$
\begin{prop}\label{prop:attainment}
Suppose that $c : \X \times \X \to \R^+$ is a bounded and continuous function, then for any $a>0$, the function $F_a$ defined in \eqref{eq:Fa} is bounded from below and attains its infimum.
\end{prop}
\proof
The cost being bounded, it is of course exponentially integrable. So according to Proposition \ref{prop:bounded-below}, the function $F_a$ is bounded below and its level sets are precompact for the weak topology.
To show that $F_a$ attains its infimum, it is enough to prove that it is lower semicontinuous (l.s.c.) with respect to the usual weak topology on $\mathcal{P}(\X).$ As it is well known, the function $\nu \mapsto H(\nu|\mu)$ is l.s.c., so it is enough to show that $\nu \mapsto \mathcal{T}_c(\nu,\mu)$ is continuous. To that end, we  take a sequence of probability measures $(\nu_n)_{n \in \N}$ converging weakly to some $\nu$, and  show that $\mathcal{T}_c(\nu_n,\mu)$ converges to $\mathcal{T}_c(\nu,\mu)$ as $n \to \infty.$
For all $n \in \N$, there exists an optimal coupling $\pi_n$ such that $\mathcal{T}_c(\nu_n,\mu) = \iint c(x,y)\,\pi_n(dxdy)$ (see \textit{e.g.} \cite[Theorem 4.1]{Vil09}). Since $\nu_n$ is a converging sequence, the sequence $\pi_n$ is tight (see \textit{e.g.} \cite[Lemma 4.4]{Vil09}) and therefore, according to Prokhorov Theorem, one can extract a subsequence $\nu_{n'}$  which converges to some coupling $\pi$ between $\nu$ and $\mu$. According to \cite[Theorem 5.20]{Vil09}, this coupling $\pi$ is also optimal, namely $\mathcal{T}_c(\nu,\mu) = \iint c(x,y)\,\pi(dxdy).$ Now, since $c$ is bounded continuous, it follows from the very definition of weak convergence, that
\[
\mathcal{T}_c(\nu_{n'},\mu) = \iint c(x,y)\,\pi_{n'}(dxdy) \to \iint c(x,y)\,\pi(dxdy)= \mathcal{T}_c(\nu,\mu),
\]
as $n' \to \infty$. Since,  by the same reasoning, any subsequence of $\nu_n$  has a subsequence $\pi_{n''}$ such that $\mathcal{T}_c(\nu_{n''},\mu)  \to  \mathcal{T}_c(\nu,\mu)$, the proof is complete.
\endproof

Now let us see how to drop the boundedness assumption on the cost.
When the cost function is of power type (as in Lemma \ref{lem:power-type}), one has the following first result.
\begin{prop}\label{prop:attainment2}
Suppose that $c:\X\times \X \to \R^+$ is a power type cost function. If $\mu$ satisfies the following strong integrability condition :
\begin{equation}\label{eq:strong-int}
\iint e^{\delta c(x,y)}\,\mu(dx)\mu(dy) <\infty,\qquad \forall \delta >0,
\end{equation}
then for all $a >0$, the function $F_a$ is bounded from below on $\mathcal{P}_\mu(\X)$ and attains its minimum.
\end{prop}
\proof
According to Proposition \ref{prop:bounded-below}, the function $F_a$ is bounded below on $\mathcal{P}_\mu(\X)$ and for all $r\in \R$, $\{F_a \leq r\}$ is precompact for the weak topology.
To conclude we need to prove  that $F_a$ is l.s.c on (say) $E:=\{F_a \leq \inf F_a + 1\}$. Since $\nu \mapsto H(\nu|\mu)$ is l.s.c on $\mathcal{P}(\X)$, it is enough to show that
$\nu \mapsto \mathcal{T}_c(\nu,\mu)$ is continuous on $E.$

\textit{First case.} First let us treat the particular case where $c(x,y)=d^p(x,y)$, $x,y\in \X$, for some $p\geq 1$ (in other words, $\mathcal{T}_c$ is the Wasserstein distance $W_p^p$). Let $(\nu_n)_{n\geq1}$ be a sequence of elements of $E$ converging to $\nu \in E$. According to \cite{Vil09}, $\mathcal{T}_c(\nu_n,\mu) \to \mathcal{T}_c(\nu,\mu)$ if and only if
\[
\limsup_{r \to \infty} \sup_{n\geq 1}\int d^p(x_o,x)\mathbf{1}_{d(x_o,x) \geq r}\,\nu_n(dx) = 0,
\]
for some (and thus all) $x_o \in \X.$
To ensure this uniform integrability condition, it is sufficient (and in fact necessary due to the de La Vallée Poussin Theorem) to prove that
\begin{equation}\label{eq:UI}
\sup_{n\geq 1} \int f(d^p(x_o,x))\,\nu_n(dx) <\infty,
\end{equation}
for some $f:\R^+ \to \R^+$ such that $f(t)/t \to \infty.$
Consider the functions
\[
\Lambda (s) = \log \int e^{s d^p(x_o,x)}\,\mu(dx),\qquad s\in \R,
\]
and
\[
\Lambda^*(t) = \sup_{s \in \R}\{st-\Lambda(s)\},\qquad t \in \R.
\]
According to \textit{e.g.} \cite[Lemma 2.2.20]{DZ98}, $\Lambda^*(t)/t \to +\infty$, as $t \to +\infty$ and,  by  \cite[Lemma 5.1.14]{DZ98},
\[
\int e^{u \Lambda^*(d^p(x_o,x))}\, \mu(dx) \leq \frac{2}{1-u},\qquad \forall u \in (0,1).
\]
According to \eqref{eq:level-sets}, there exists $h>0$ such that $E\subset \{\nu \in \mathcal{P}(\X) : H(\nu|\mu) \leq h\}$. So, using Young's inequality \eqref{eq:Young} as in Proposition \ref{prop:bounded-below}, one easily gets that for all $u\in (0,1)$,
\[
\int \Lambda^*(d^p(x_o,x))\,\nu_n(dx) \leq \frac{1}{u} H(\nu_n|\mu) + \frac{e^{-1}}{u} \int e^{u\Lambda^*(d^p(x_o,x))}\,\mu(dx) \leq \frac{h}{u} + \frac{e^{-1}}{u} \frac{2}{1-u}.
\]
This shows \eqref{eq:UI} with $f = \Lambda^*$ and completes the proof of the first case.

\textit{Second case.}  Now we treat the general case. According to Lemma \ref{lem:power-type}, the function $\tilde{d} = c^{1/p_o}$ is a distance on $\X$, where $p_o$ denotes the exponent of $c$ and is defined in the lemma. So we have $c(x,y) = \tilde{d}^p(x,y)$, for all $x,y \in \X$ and we are back to the first case (since according to Lemma \ref{lem:power-type} $(\X,\tilde{d})$ is still a Polish space). Details are left to the reader.
\endproof

The strong integrability condition \eqref{eq:strong-int} is a bit too demanding for our purpose.
For instance, in the particular case of  the functional
\[
F_a(\nu) = aH(\nu|\mu) - W_2^2(\nu,\mu),\qquad \nu \in \mathcal{P}_\mu(\X),
\]
existence of a minimizer is granted by the previous result if  $\iint e^{\delta d^2(x,y)}\,\mu(dx)\mu(dy)<\infty$ for all $\delta>0$. This condition is for instance not satisfied by the standard Gaussian probability measure when $\X=\R^d$, yet a minimizer (the standard Gaussian probability measure itself) in that case exists, if $a\geq 2$.

Fortunately, this integrability condition can be relaxed as shown by the following result.
\begin{thm}\label{thm:attainment3}
Let $c:\X \times \X \to \R^+$ be a power type cost function and $p_o$ be its exponent. Suppose moreover that $\mu$ satisfies the concentration property \eqref{eq:concentration} for some $a'>0$ and $r_o\geq0$. Then, for all $a>a'$, the functional $F_a =aH(\nu|\mu) - \mathcal{T}_\cost(\nu,\mu)$ admits at least one minimizer on $\mathcal{P}_\mu(\X).$
\end{thm}

So, with the notations introduced in Remark \ref{rem:comparison}, for all $a>a_o=a_o'$ the functional $F_a$ is lower bounded and attains its minimum. For all $a<a_o=a_o'$ the functional $F_a$ is not bounded from below.

The proof of Theorem \ref{thm:attainment3} is sensibly more sophisticated than those of Propositions \ref{prop:attainment} and \ref{prop:attainment2} (and actually relies on them). It is postponed to Section \ref{Sec:attainment}. Let us mention, that in order to stay at the most elementary level, we prefer to use in the applications considered in Section \ref{Sec:Applications} a simple truncation of the cost technique which will enable us to use Proposition \ref{prop:attainment} instead of Theorem \ref{thm:attainment3}.

\subsection{Characterization of the minimizers}
According to \textit{e.g.} \cite[Theorem 5.10]{Vil09}, we have the following Kantorovich duality formula
\begin{equation}\label{duality}
\mathcal{T}_\cost(\nu,\mu) = \sup \left\{ \int \psi\,d\nu + \int \varphi\,d\mu \right\},
\end{equation}
where the supremum runs over all functions $\psi$ and $\varphi$ such that $\psi \in L^1(\nu), \varphi \in L^1(\mu)$ and
\[
\psi(x)+\varphi(y)\leq \cost(x,y),\qquad \forall (x,y)\in\X^2.
\]
Our first task is to show that there exists an optimal pair $(\overline{\psi},\overline{\varphi})$ for any $\nu \in \mathcal{P}_\mu(\X).$
\begin{prop}\label{prop:optimal-pair}
Assume that $\mu$ satisfies the exponential integrability condition \eqref{eq:expmoment} for some $\delta>0$. Then for all $\nu \in \mathcal{P}_\mu(\X)$, there exist $\overline{\psi} \in L^1(\nu)$ and $\overline{\varphi} \in L^1(\mu)$ related together by the
$\cost$-conjugation:
\begin{equation}\label{eq:c-trans}
\overline{\psi}(x)=\inf_{y\in\X}\{\cost(x,y)-\overline{\varphi}(y)\}, x \in \X\qquad \overline{\varphi}(y)=\inf_{x\in\X}\{\cost(x,y)-\overline{\psi}(x)\}, y\in \X.
\end{equation}
such that
\[
\mathcal{T}_c(\nu,\mu) = \int \overline{\psi}(x)\,\nu(dx) + \int \overline{\varphi}(y)\,\mu(dy).
\]
\end{prop}
We recall that a function $\overline{\psi}$ as in the proposition above, is called a \emph{Kantorovich potential} for the transport of $\nu$ on $\mu$.
(Note that $\overline{\varphi}$ is determined by $\overline{\psi}$ according to \eqref{eq:c-trans}.)

\proof
According to Proposition \ref{prop:bounded-below}, $\iint c(x,y)\,\nu(dx)\mu(dy)$ is finite for any $\nu \in \mathcal{P}_\mu(\X)$. According to \cite[Theorem 6.1.5 and Remark 6.1.6]{AmbGigSav-08}, this is enough to ensure the existence of the desired optimal pair.
\endproof

We are now ready to state our characterization of minimizers of $F_a$.
\begin{thm}\label{thm:CharactMinimizer}
Let $\mu$ be a probability measure on $\X$ satisfying the integrability condition  \eqref{eq:expmoment} for some $\delta>0.$  Assume that for some value of $a>0$, the function $F_\const$ is bounded from below and reaches its infimum at some point $\underline{\nu}\in \mathcal{P}_\mu(\X)$ not equal to $\mu$. Then the density of $\underline{\nu}$ with respect to $\mu $ satisfies the following equation $\underline{\nu}$-almost everywhere in $\X$
\[
\overline{\lambda} \log\left(\frac{d\underline{\nu}}{d\mu}\right) = \overline{\psi}+C,\qquad \text{with}\qquad \overline{\lambda} = \frac{a\alpha'(aH(\underline{\nu}|\mu))}{\beta'(\mathcal{T}_c(\underline{\nu}|\mu))}
\]
where $C\in \R$ is a renormalizing constant and $\overline{\psi}$ is any Kantorovich potential for the transport of $\underline{\nu}$ on $\mu$.
\end{thm}

First let us state  two useful  lemmas  which follow ideas in \cite{MauRoudSantam} further developed in \cite{Santambook}  ($\mathcal{C}_b(\X)$ denotes here  the set of bounded continuous functions on $\X.$)
\begin{lem}\label{dif-H}
Let $\underline{\nu}$ be absolutely continuous with respect to $\mu$ and $H(\underline{\nu}|\mu)<\infty.$
Let $f \in\mathcal{C}_b(\X)$ be such that $\int f\,d\underline{\nu}=0$ and define $h(\varepsilon) = H(\nu_{\varepsilon} |\mu)$ with $\nu_\varepsilon = (1+\varepsilon f)\underline{\nu}$.  Then
\[
h'(0) = \int \log \left(\frac{d\underline{\nu}}{d\mu}\right)\,f\,d\underline{\nu}
\]
\end{lem}
\proof
The proof is left to the reader.
\endproof

\begin{lem}\label{dif-OptCost}
Let $\mu,\underline{\nu} \in \mathcal{P}(\X)$ and suppose that $\overline{\varphi} \in L_1(\underline{\nu}),\overline{\psi} \in L_1(\mu)$ is a couple of functions such that $\overline{\psi}(x) + \overline{\varphi}(y) \leq \cost(x,y)$ for all $(x,y) \in \X^2$ and
\begin{equation}\label{potential}
\mathcal{T}_\cost(\underline{\nu},\mu) = \int \overline{\psi}\,d\underline{\nu} + \int \overline{\varphi}\,d\mu.
\end{equation}
Then, for all $\nu \in \mathcal{P}(\X)$ such that $\overline{\psi} \in L_1(\nu)$, it holds
\[
\mathcal{T}_\cost(\nu,\mu) \geq \mathcal{T}_\cost(\underline{\nu},\mu) + \int \overline{\psi}\,d(\nu-\underline{\nu}).
\]
\end{lem}
\proof
 By Kantorovich duality $\mathcal{T}_\cost(\nu,\mu) \geq \int \psi\,d\nu + \int \varphi\,d\mu$, for all couple of integrable functions such that $\psi(x)+\varphi(y) \leq \cost(x,y)$ for all $x,y\in \X.$ Therefore, taking the couple $(\overline{\psi},\overline{\varphi})$, with $\overline{\varphi}(y) = \inf_{x}\{-\overline{\psi}(x)+c(x,y)\}$, $y\in \X$  immediately  yields
\[
\mathcal{T}_\cost(\nu,\mu)  \geq \int \overline{\psi}(x)\,\nu(dx) + \int \overline{\varphi}(x)\,\mu(dx)= \mathcal{T}_\cost(\underline{\nu},\mu) + \int \overline{\psi}\,d(\nu-\underline{\nu}).
\]
\endproof

\proof[Proof of Theorem \ref{thm:CharactMinimizer}]
It holds, for all $\nu\in \mathcal{P}(\X)$
\[
F_\const(\nu) -F_\const(\underline{\nu})= \alpha\left(\const H(\nu|\mu)\right)-\alpha\left(\const H(\underline{\nu}|\mu)\right)-\beta\left(\mathcal{T}_\cost(\nu,\mu)\right)+\beta\left(\mathcal{T}_\cost(\underline{\nu},\mu)\right).
\]
Since $\mu$ satisfies the integrability condition \eqref{eq:expmoment} for some $\delta>0$ and $\underline{\nu} \in \mathcal{P}_\mu(\X)$, Proposition \ref{prop:optimal-pair} shows that there
exists an optimal pair $(\overline{\psi},\overline{\varphi})$ for the transport of $\underline{\nu}$ on $\mu$. Since $\beta$ is non-decreasing, using Lemma \ref{dif-OptCost} we get
\[
-\beta\left(\mathcal{T}_\cost(\nu,\mu)\right)\leq -\beta\left(\mathcal{T}_\cost(\underline{\nu},\mu) +\int \overline{\psi}\,d(\nu-\underline{\nu})\right).
\]
Therefore, taking $\nu = \nu_\varepsilon = (1+\varepsilon f)\,\underline{\nu}$ for functions $f \in \mathcal{C}_b(\X)$ such that $\int f\,d\underline{\nu}=0$ and letting $\varepsilon \to 0$ yields, according to Lemma \ref{dif-H},
\begin{align*}
\limsup_{\varepsilon\to 0+} \frac{F_\const(\nu_\varepsilon) -F_\const(\underline{\nu})}{\varepsilon} &\leq \lim_{\varepsilon\to 0+} \frac{\alpha\left(\const H(\nu_\varepsilon|\mu)\right)-\alpha\left(\const H(\underline{\nu}|\mu)\right)}{\varepsilon} \\
&+
\lim_{\varepsilon\to 0+} \frac{\beta\left(\mathcal{T}_\cost(\underline{\nu},\mu)\right)-\beta\left(\mathcal{T}_\cost(\underline{\nu},\mu) +\varepsilon \int \overline{\psi}\,f d\underline{\nu}\right)}{\varepsilon}\\
&\,= \const \alpha'\left(\const H(\underline{\nu}|\mu)\right)\int \log \left(\frac{d\underline{\nu}}{d\mu}\right)\,f\,d\underline{\nu}-\beta'\left(\mathcal{T}_\cost(\underline{\nu},\mu)\right)\int \overline{\psi}f\,d\underline{\nu}.
\end{align*}
Now, since $F_\const$ reaches its infimum at $\underline{\nu}$, the left hand side is non-negative, and so
\[
\int \left(\const \alpha'\left(\const H(\underline{\nu}|\mu)\right)\log \left(\frac{d\underline{\nu}}{d\mu}\right)- \beta'\left(\mathcal{T}_\cost(\underline{\nu},\mu)\right)\overline{\psi}\right)\,f\,d\underline{\nu} \geq 0,
\]
for all $f\in \mathcal{C}_b(\X)$ such that $\int f\,d\underline{\nu} =0.$ Changing $f$ into $-f$, one concludes that there is in fact equality. We conclude applying the lemma below.
\endproof

\begin{lem}\label{nul}
Let $\nu \in \mathcal{P}(\X)$ and suppose that $g \in L_1(\nu)$ is such that $\int fg\,d\nu = 0$ for all $f \in \mathcal{C}_b(\X)$ such that $\int f\,d\nu =0.$ Then $g$ is constant $\nu$-almost surely.
\end{lem}
\proof
First observe that   $\int \left(g-\int g\,d\nu\right)f\,d\nu =\int \left(g-\int g\,d\nu\right) \left(f-\int f\,d\nu\right)\,d\nu =0$ for all $f \in \mathcal{C}_b(\X)$, hence the signed Borel measure $ \left(g-\int g\,d\nu\right)\cdot d\nu$  is null. In particular,  for every  $t>0$ it holds that
\[
\int \left(g-\int g\,d\nu\right)\mathbf{1}_A\,d\nu =0
\]
for $A = \{g - \int g\,d\mu >t\}$ and  $A = \{g - \int g\,d\mu <-t\}$ which yields $g = \int g\,d\nu$, $\nu$-almost surely.
\endproof

Now we consider the special case where $\X=\R^d$ and $c(x,y) = |x-y|^2$, $x,y \in \R^d$, where $|\,\cdot\,|$ denotes the standard Euclidean norm. In this situation, minimizers of $F_a$ are solutions of a certain Monge-Amp\`ere equation, as shown in the following result.
\begin{cor}\label{cor:MA} Let $\mu$ be a probability measure on $\R^d$ absolutely continuous with respect to Lebesgue measure and satisfying the integrability condition  \eqref{eq:expmoment} for some $\delta>0.$  Assume that for some value of $a>0$, the function $F_\const$ is bounded from below and reaches its infimum at some point $\underline{\nu}\in \mathcal{P}_\mu(\X)$ not equal to $\mu$. Then
$\underline{\nu}$ admits a density with respect to $\mu$ of the form $\underline{\nu}(dx) = e^{-V(x)}\,\mu(dx)$   with $V$ such that $x\mapsto V(x) + |x|^2/\overline{\lambda}$ is convex on $\R^d$ and  $\overline{\lambda} = \frac{a\alpha'(aH(\underline{\nu}|\mu))}{\beta'(\mathcal{T}_c(\underline{\nu}|\mu))}$.
As such, $V$ is differentiable and admits a Hessian in the sense of Aleksandrov, Lebesgue almost everywhere.
Moreover denoting by $T$ the Brenier map sending $\underline{\nu}$ to $\mu$, then for $\underline{\nu}$-almost all $x$, it holds
\begin{equation}\label{eq:preMA}
\overline{\lambda} \nabla V(x) = 2 (T(x)-x).
\end{equation}
Finally, $V$ is a solution to the following Monge-Ampère type equation: for $\underline{\nu}$-almost all $x \in \R^d$,
\begin{equation}\label{eq:MA}
h\left(x+\frac{\overline{\lambda}}{2} \nabla V (x)\right)\det\left(I_d+\frac{\overline{\lambda}}{2}\nabla^{2}V(x)\right)=e^{-V(x)}h(x),
\end{equation}
where $h$ is the density of $\mu$ with respect to Lebesgue measure and $I_d$ is the identity matrix in $\R^d$.
\end{cor}

\proof[Proof of Corollary \ref{cor:MA}]
According to Theorem \ref{thm:CharactMinimizer}, it holds for $\underline{\nu}$-almost every $x\in \R^d$,
$\overline{\lambda} \log\left(\frac{d\underline{\nu}}{d\mu}\right)(x) = \overline{\psi}(x)+C$, for some $C\in \R$, where $\overline{\psi}$ is a Kantorovich potential for the transport of $
\underline{\nu}$ on $\mu.$ Modifying the density of $\underline{\nu}$ on a negligible set, one can assume that this equality holds everywhere. Setting $V(x) = -(\overline{\psi}(x) + C)/\overline{\lambda}$, the equation reads $\underline{\nu}(dx) = e^{-V(x)}\,\mu(dx)$.
According to Proposition \ref{prop:optimal-pair}, there is an adjoint function $\overline{\varphi}$ such that
\[
\overline{\psi}(x) = \inf_{y\in \R^d}\{-\overline{\varphi}(y)+|x-y|^2\}
 = |x|^2 -\sup_{y\in \R^d}\{\overline{\varphi}(y)-|y|^2 + 2 x\cdot y\}.
\]
Being a supremum of linear functions, the last function is convex, which shows that $x \mapsto V(x) + |x|^2/\overline{\lambda}$ is convex on $\R^d$. According to Aleksandrov Theorem (see \textit{e.g.} Evans-Gariepy \cite{EvanGari-92}), it follows that for all $x \in \R^d$ outside a set of Lebegue measure $0$, the function $V$ is differentiable at $x$ and there is a symmetric matrix denoted $\nabla^2 V(x)$ such that
\[
V(x+h) = V(x) + \nabla V(x)\cdot h + \frac{1}{2} \nabla^2V(x)h\cdot h+ o(|h|^2),
\]
as $h \to 0.$ According to \cite[Theorem 1.2]{GMC96}, the Brenier map $T$ transporting $\underline{\nu}$ on $\mu$ is related to $\overline{\psi}$ by $2(x-T(x)) = \nabla \overline{\psi}(x)$ for $\underline{\nu}$-almost all $x \in \R^d$, which gives \eqref{eq:preMA}.
According to the change of variable formula (see \cite[Theorem 4.4 and Remark 4.5]{MC97}), for $\underline{\nu}$-almost all $x \in \R^d$,
\begin{equation}\label{ChgVar}
h(T(x))\det(DT(x))=e^{-V(x)},
\end{equation}
where $DT$, the differential of $T$, is expressed at every point where it is well defined, by
$DT(x) = I_d + \frac{\overline{\lambda}}{2}\nabla^2V(x)$. This completes the proof.
\endproof

\section{Applications}\label{Sec:Applications}
In this section, we use our characterization of minimizers to recover different known implications between classical functional inequalities.

Recall that if $\nu$ is absolutely continuous with respect to $\mu$, the Fisher information $I(\nu|\mu)$ of $\nu$ with respect to $\mu$ is by definition
\[
I(\nu | \mu) = \int \left|\nabla^{+} \log \left( \frac{d\nu}{d\mu}\right)\right|^2\, d\nu.
\]
(Recall the definition of the local slope $|\nabla ^+ g|$ for a function $g:\X \to \R$ given in \eqref{eq:local-slope}.)
\subsection{A new proof of Otto-Villani Theorem}
The following result is due to Otto and Villani \cite{OV00}.
\begin{thm}\label{thm:OV}
Suppose that $\mu$ satisfies the logarithmic Sobolev inequality
\[
H(\nu|\mu) \leq D I(\nu|\mu), \qquad \forall \nu \in \mathcal{P}(\X),
\]
then $\mu$ satisfies $\T_2(4D).$
\end{thm}
\proof
\textit{First case.} First we treat the case where the distance $d$ is bounded on $\X^2.$
Let us consider
\begin{equation}\label{Gc}
F_\const(\nu) = \sqrt{\const H(\nu|\mu)} - W_2(\nu,\mu),
\end{equation}
for $\const>4D$, which is well defined on $\mathcal{P}(\X).$
According to Proposition \ref{prop:attainment}, the function $F_a$ reaches its infimum at some point $\underline{\nu} \in \mathcal{P}(\X).$ Our goal is to prove that $\underline{\nu} = \mu.$ Let us assume, by contradiction, that $\underline{\nu} \neq \mu.$ According to Theorem \ref{thm:CharactMinimizer},
\[
\frac{\sqrt{\const}h(x)}{\sqrt{ H(\underline{\nu}|\mu)}}=\frac{\overline{\psi}(x)}{W_2(\underline{\nu},\mu)}+C,
\]
for $\underline{\nu}$-almost every $x\in\X$,
where  $h(x)=\log\left(\frac{d\underline{\nu}}{d\mu}\right)(x).$
Modifying $\frac{d\underline{\nu}}{d\mu}$ on a set of $\underline{\nu}$-null measure, we can assume that the preceding equality holds true for all $x$ in $\X$.
In particular, taking the local slope, it holds
\[
|\nabla^+ h |(x) = \frac{\sqrt{H(\underline{\nu}|\mu)}}{\sqrt{\const}W_2(\underline{\nu},\mu)} |\nabla^+\overline{\psi}|(x),\qquad \forall x\in\X.
\]
Therefore, we have
\[
\int |\nabla^+ h |^2(x)\,\underline{\nu}(dx) =\frac{H(\underline{\nu}|\mu)}{\const W^2_2(\underline{\nu},\mu)} \int |\nabla^+\overline{\psi}|^2(x)\,\underline{\nu}(dx)
\]
Applying the following lemma (whose proof is given below),
\begin{lem}\label{lem1} Let $\nu,\mu\in \Pp(\X)$ have finite second moments and  $\psi$ be a Kantorovich potential for the transport of $\nu$ on $\mu$ (with respect to the cost $c(x,y) = d^2(x,y)$). Then,
\[
\int |\nabla^+\psi|^2(x)\,\nu(dx)\leq 4 W^2_2(\nu,\mu).
\]
\end{lem}
\noindent it follows that
\[
I(\underline{\nu}|\mu) =\int |\nabla^+h|^2(x)\,\underline{\nu}(dx) \leq \frac{4 H(\underline{\nu}|\mu)}{\const}.
\]
Now since $\mu$ satisfies the logarithmic Sobolev inequality
\[
H(\underline{\nu}|\mu) \leq D I(\underline{\nu}|\mu)
\]
one concludes that $1 \leq (4D)/\const$, which is impossible. Hence $\underline{\nu} =\mu.$
Therefore $F_\const(\nu) \geq F_\const(\mu) =0$ for all $\nu$ which proves that $\mu$ satisfies $\T_2(\const)$ for all $\const>4D$.  Letting $\const$ go to $4D$, this ends the proof of the first case.

\noindent \textit{Second case.} Now we treat the case where the distance $d$ is unbounded on $\X^2.$ For all positive integer $n$, let us consider the distance $d_n$ defined by
\[
d_n(x,y) = d(x,y) \wedge n,\qquad \forall (x,y) \in \X^2.
\]
Observe that, for any function $f:\X \to \R$, it holds
\[
\limsup_{y \to x} \frac{[f(y)-f(x)]_+}{d_n(x,y)} = \limsup_{y \to x} \frac{[f(y)-f(x)]_+}{d(x,y)}.
\]
Therefore, one immediately concludes that $\mu$ satisfies the logarithmic Sobolev inequality with constant $D$ on the metric space $(\X,d_n)$. Applying the first case, we conclude that $\mu$ satisfies the inequality
\[
W_{2,n}^2(\nu,\mu) \leq 4 D H(\nu|\mu),\qquad \forall \nu \in \mathcal{P}(X),
\]
where $W_{2,n}^2(\nu,\mu) = \inf_{\pi} \iint d_n^2(x,y)\,\pi(dxdy)$,  the infimum running over all couplings between $\nu$ and $\mu.$ To complete the proof, it suffices to show that
\begin{equation}\label{eq:liminf-utile}
W_{2}^2(\nu,\mu) \leq \liminf_{n \to \infty}W_{2,n}^2(\nu,\mu).
\end{equation}
For all $n\geq 1$, consider an optimal coupling $\pi_n$ of  $\nu$ and $\mu$ for the cost $d_n^2$. The sequence $\pi_n$ is tight, and so some henceforth equally denoted subsequence   converges weakly to some coupling $\pi$ of $\nu$ and $\mu$.
Fix $m \geq 1$, then for all $n \geq m$, it holds
\[
\iint d_m^2(x,y)\,\pi_n(dxdy) \leq \iint d_n^2(x,y)\,\pi_n(dxdy).
\]
Therefore taking the $\liminf$ when $n \to \infty$ and using the weak convergence of the sequence $\pi_n$, one gets
\[
\iint d_m^2(x,y) \,\pi(dxdy) \leq \liminf_{n \to \infty}W_{2,n}^2(\nu,\mu).
\]
Using monotone convergence theorem as $m$ goes to $\infty$, one then concludes that
\[
W_2^2(\nu,\mu) \leq \iint d^2(x,y)\,\pi(dxdy)  \leq \liminf_{n \to \infty}W_{2,n}^2(\nu,\mu),
\]
which completes the proof of the second case.
\endproof
\proof [Proof of Lemma \ref{lem1}]
Let $\psi$ be a Kantorovich potential between $\nu$ and $\mu$ (whose existence is given by \textit{e.g.} \cite[Theorem 5.10]{Vil09}); by definition
\begin{equation}\label{eq:proofLem1}
\psi(x)=\inf_{y\in\X}\{d^2(x,y)-\varphi(y)\}.
\end{equation}
with $\varphi (y) = \inf_{x \in \X}\{-\psi(x) + d^2(x,y)\}$, $y\in \X$ (in the terminology of optimal transport, the function $\psi$ is said $d^2$-concave).

Recall that the $d^2$-subdifferential of $\psi$ at a point $x \in \X$ is the (possibly empty) set denoted by $\partial_{d^2}\psi(x)$ of points $\bar{y}$ realizing the infimum in \eqref{eq:proofLem1}. If $\bar{y} \in \partial_c\psi(x)$, we get
\[
\psi(z) - \psi(x) \leq d^2(z,\bar{y}) - d^2(x,\bar{y}),\qquad \forall z \in \X.
\]
(Note that this inequality is an equivalent definition of $\bar{y} \in \partial_c \psi(x)$.)
Hence, for all $\overline{y}$ in $\partial_{d^2}\psi(x)$ (when non empty), we have
\begin{align*}
\psi(z)-\psi(x)&\leq \left(d(z,\overline{y})-d(x,\overline{y})\right)\left(d(z,\overline{y})+d(x,\overline{y})\right)\\
&\leq d(z, x)\left(d(z,\overline{y})+d(x,\overline{y})\right).
\end{align*}
Therefore, for all $x \in \X$,
\[
|\nabla^+\psi|(x)\leq 2\inf_{\overline{y}\in\partial_{d^2}\psi(x)}d(x,\overline{y})
\]
(with the convention $\inf_{\emptyset} = +\infty$)
and thus
\[
\int |\nabla^+\psi|^2(x)\,\nu(dx)\leq 4\int \inf_{\overline{y}\in\partial_{d^2}\psi(x)}d^2(x,\overline{y})\,\nu(dx).
\]
Let us denote by
\[
\partial_c \psi = \{(x,\bar{y}) : x\in \X, \bar{y} \in \partial_c\psi(x)\}\subset \X^2.
\]
Let $\pi$ be an optimal coupling between $\nu$ and $\mu$ (whose existence is given by \textit{e.g.} \cite[Theorem 4.1]{Vil09}), then $\pi(\partial_{d^2}\psi)=1$ (which implies in particular that $\partial_c \psi(x) \neq \emptyset$ for $\nu$-almost all $x$).
Let us briefly justify this well known property. By optimality of $\pi$ and of $(\psi,\varphi)$, it holds
\[
0 = \iint_{\X^2} d^2(x,y)\,-(\psi(x) + \varphi(y))\,\pi(dxdy),
\]
which, since the integrand is non-negative, implies that $\psi(x) = -\varphi(y) + d^2(x,y)$ for $\pi$-almost all $(x,y) \in \X^2$, meaning  exactly  that $\pi(\partial_{d^2}\psi)=1.$
Using this fact, it thus holds that
\[
W^2_2(\nu,\mu)=\iint d^2(x,y)\,\pi(dxdy)\geq \int \inf_{\overline{y}\in\partial_{d^2}\psi(x)}d^2(x,\overline{y})\,\nu(dx)
\]
which ends the proof.
\endproof

\subsection{A variant involving a restricted logarithmic Sobolev inequality}
The following theorem is a variant of a result obtained by the second author in a paper with Roberto and Samson (see \cite{GRS11b} for the case of the Euclidean case and \cite{GRS13,GRS14} for the metric space case).

For $\lambda_o>0$, let $\mathcal{F}_{\lambda_o}(\X)$ denote the class of $\lambda d^2$-concave functions with $\lambda \in (0,\lambda_o)$, \textit{i.e.} the set of functions $f:\X \to \er$ for which there exists $g:\X \to \er$  and $\lambda \in (0,\lambda_o)$ such that
\[
f(x) = \inf_{y \in \X} \{ -g(y) + \lambda d^2(x,y)\},\qquad x \in \X.
\]
Let us remark that if $\X = \er^d$ is equipped with the usual Euclidean norm $|\,\cdot\,|$, then a function $f$ belongs to $\mathcal{F}_{\lambda_o}(\R^d)$ if and only if the function $x\mapsto f(x) - \lambda |x|^2$, $x \in \R^d$, is concave, which translates into the following semi-concavity property
\begin{equation}\label{eq:semiconv}
f((1-t) x + t y) \geq (1-t)f(x) + t f(y) - \lambda t(1-t)|x-y|^2, \qquad \forall x,y \in \R^d,\qquad \forall t\in [0,1].
\end{equation}

\begin{thm}\label{thm:variant-GRS}
Suppose that for some constants $\lambda_o, D >0$, $\mu$ satisfies the following restricted version of the logarithmic Sobolev inequality
\[
H(\nu|\mu) \leq D I(\nu|\mu),
\]
for all $\nu$ of the form $\nu(dx) = e^{f(x)}\,\mu(dx)$ with $f$ belonging to the class $\mathcal{F}_{\lambda_o}(\X)$ defined above. Then $\mu$ satisfies $\T_2(\max(4D,1/\lambda_o)).$
\end{thm}

In comparison, it was shown in \cite{GRS11b} that a probability measure $\mu$ on a Euclidean space which satisfies the logarithmic Sobolev inequality $H(\nu|\mu) \leq D I(\nu|\mu)$, for all $\nu(dx)=e^{-f(x)}\,\mu(dx)$ with  $f \in \mathcal{F}_{\lambda_o}(\X)$ also satisfies $\T_2(C)$, for some $C$ depending on $D$ and $\lambda_o$. Moreover, the converse is also true : if $\mu$ satisfies $\T_2(C)$ it also satisfies this restricted logarithmic Sobolev inequality for some $D$ and $\lambda_o$. See \cite{GRS13,GRS14} for extensions to general geodesic spaces. We do not know if there is also equivalence in Theorem \ref{thm:variant-GRS}.

\proof
\textit{First case.} Assume that the distance $d$ is bounded. Take $a>\max (4D;\lambda_o)$ and consider again the function $F_a$ defined by \eqref{Gc}. Suppose that $F_a$ reaches its minimum at some point $\underline{\nu} \neq \mu$. Reasoning as in the proof of Theorem \ref{thm:OV}, one concludes that the density $\rho$ of $\underline{\nu}$ satisfies (possibly after a modification on a set of $\underline{\nu}$ measure $0$)
\[
\log \rho(x) =  \lambda  \overline{\psi}(x),\qquad \forall x\in \X,
\]
where $\overline{\psi} \in \mathcal{F}_1(\X)$ and $\lambda=\frac{\sqrt{H(\underline{\nu}|\mu)}}{\sqrt{a}W_2(\underline{\nu},\mu)}$.
Since $\underline{\nu}\neq \mu$, $F_a(\underline{\nu}) \leq F_a(\mu) = 0$ and so $\sqrt{a H(\underline{\nu}|\mu)} \leq W_2(\underline{\nu},\mu).$ Therefore, $\lambda \leq 1/a < \lambda_o$ and so $\log \rho \in \mathcal{F}_{\lambda_o}(\X).$
Since $\mu$ satisfies the logarithmic Sobolev  inequality for such densities, reasoning as in the proof of Theorem \ref{thm:OV}, we conclude exactly as before that $\underline{\nu}$ can not be distinct from $\mu$. Therefore $\mu$ satisfies $\T_2(a)$. Letting $a$ go to $\max (4D;\lambda_o)$ completes the proof of the first case.\\
\textit{Second case.} Now we assume that $d$ is unbounded and we use the same truncation trick as in the proof of Theorem \ref{thm:OV}. In order to reason exactly as before, all what we need to check is that if $f,g$  are functions such that
\[
f(x) = \inf_{y \in \X}\{ -g(y) + \lambda (d(x,y)\wedge M)^2\}, \quad x \in \X
\]
for fixed $M>0$ then $f$ belongs to $\mathcal{F}_\lambda (\X).$ Define $h(y)=-\inf\{ -f(x)+\lambda (d(x,y)\wedge M)^2\}$, $y\in \X.$ Then it is easily checked that $f(x) =  \inf_{y \in \X}\{ -h(y) + \lambda (d(x,y)\wedge M)^2\}$,  $x \in \X.$

Let us show that the function $h$ is bounded and satisfies
\begin{equation}\label{eq:osch}
h(y)-h(x) \leq \lambda M^2, \qquad \forall x,y \in \X.
\end{equation}
By definition of $h$, for all $x \in \X$, $h(x) \geq \inf (-f) = - \sup f$, and so $\inf h \geq -\sup f.$ On the other hand, $f(y) \leq  \inf (-h) + \lambda M^2 = - \sup (h)+\lambda M^2$, and so $\sup f \leq -\sup h + \lambda M^2.$ We conclude from this that $\sup h \leq \inf h +\lambda M^2$, which amounts to \eqref{eq:osch}.

Now, let us define $\tilde{f}(x) = \inf_{y\in \X}\{ -h(y) + d^2(x,y)\}$, $x \in \X$ and let us show that $\tilde{f}=f.$ Fix a point $x \in \X$ and first observe that $f(x) \leq -h(x)$ and $\tilde{f}(x) \leq -h(x)$.
On the other hand, if $d(x,y)\geq M$, then it follows from \eqref{eq:osch} that  $-h(y) + \lambda d(x,y)^2\geq -h(x).$ Therefore,
\[
\tilde{f}(x) = \inf_{y \text{ s.t. } d(x,y)< M }\{ -h(y) + \lambda d^2(x,y)\} =  \inf_{y \text{ s.t. } d(x,y)< M }\{ -h(y) + \lambda (d(x,y)\wedge M)^2\}.
\]
Similarly, $f(x) = \inf_{y \text{ s.t. } d(x,y)< M }\{ -h(y) + \lambda (d(x,y)\wedge M)^2\}$ and so $f(x) = \tilde{f}(x).$\endproof

\subsection{A variant involving a transport-information inequality}

In the following result  we show that the transport-information inequality $W_2I$ introduced by Guillin, L\'eonard, Wu and Yao \cite{GLWY09} implies $\T_2$ (see also \cite{GLWW09} for a proof using the Hamilton-Jacobi method).
\begin{prop}\label{prop:WI}
Suppose that $\mu$ satisfies the inequality
\[
W_2^2(\nu,\mu) \leq D I(\nu|\mu),
\]
for all $\nu,$ where $I(\nu | \mu) = \int |\nabla^+ \log \left( \frac{d\nu}{d\mu}\right)|^2\, d\nu$ denotes the Fisher information.
Then $\mu$ satisfies $\T_2(2\sqrt{D}).$
\end{prop}
\proof
Here, we will only consider the case where the metric $d$ is bounded. The case of unbounded metrics is treated exactly as in the proof of Theorem \ref{thm:OV}. Take $\const >2\sqrt{D}$ and consider the functional $F_\const$ defined by
\[
F_\const(\nu) = \const H(\nu|\mu) -W_2^2(\nu,\mu),\qquad \forall \nu \in \mathcal{P}(\X).
\]
Since, the metric is bounded, the function $F_a$ is bounded from below and so according to Proposition \ref{prop:attainment} and Theorem \ref{thm:CharactMinimizer}, it reaches its infimum at some point $\underline{\nu}$ whose density $\rho$ with respect to $\mu$  satisfies the following equality
\[
\const\log \rho = \overline{\psi}+C,
\]
where $\overline{\psi}$ is a Kantorovich potential between $\underline{\nu}$ and $\mu$ and $C$ some constant.
Taking the local slope and reasoning as in the proof of Theorem \ref{thm:OV}, one concludes that
\[
4 W_2^2(\underline{\nu},\mu) = \const^2I(\underline{\nu}|\mu).
\]
By assumption, $I(\underline{\nu}|\mu) \geq \frac{1}{D} W_2^2(\underline{\nu}|\mu)$. Since $\frac{\const^2}{D} >4$, one concludes that $\underline{\nu} = \mu.$ It follows that $F_\const(\nu) \geq F_\const(\mu) = 0$ for all $\nu$ and so $\mu$ satisfies $\T_2(\const)$. Letting $\const \to 2\sqrt{D}$ completes the proof.
\endproof
\subsection{Comparison between optimal constants}
For a given probability measure $\mu$ on $(\X,d)$, we denote by $C_{LSI}(\mu)\in [0,\infty]$, $C_{T_2}(\mu)$ and $C_{W_2I}(\mu)$ the best constants (\textit{i.e.} the smallest) in the logarithmic Sobolev, in Talagrand and in the $W_2I$ inequalities for $\mu$.
Otto-Villani Theorem can be simply restated as the inequality
\[
C_{T_2}(\mu) \leq 4C_{LSI}(\mu).
\]
Combining $\LSI$ and $\T_2$, one easily sees that
\[
C_{W_2I}(\mu) \leq 4C_{LSI}^2(\mu),
\]
and according to Proposition \ref{prop:WI},
\[
C_{T_2}(\mu) \leq 2\sqrt{C_{W_2I}}(\mu).
\]
Therefore, one has the inequalities
\[
C_{T_2}(\mu) \leq 2\sqrt{C_{W_2I}(\mu)} \leq 4 C_{LSI}(\mu)
\]

Now let us assume that $\X = \R^d$ equipped with its usual Euclidean norm and that $\mu$ is absolutely continuous with respect to Lebesgue measure with a density denoted by $h:\R^d\to[0,\infty).$ We also assume, for simplicity,  that $\mu$ is \emph{compactly supported}.
Let us also denote by $\mathcal{V}_a$, the class of functions $V:\R^d\to\R\cup\{+\infty\}$ such that $x \mapsto V(x) + \frac{|x|^2}{a}$ is convex.
Consider again the functional
\[
F_a (\nu)= aH(\nu|\mu) - W_2^2(\nu,\mu),\qquad \nu \in \mathcal{P}_\mu(\X),
\]
used in Proposition \ref{prop:WI}. According to Proposition \ref{prop:attainment} (which applies since the support of $\mu$ is compact), for all $a>0$, $\mathrm{Argmin}(F_a) \neq \emptyset.$
Moreover, according to Theorem \ref{thm:CharactMinimizer} and Corollary \ref{cor:MA}, if $\underline{\nu} \in \mathrm{Argmin}(F_a)$, then it is of the form $\underline{\nu}(dx) = e^{-V(x)}\,\mu(dx)$ for some function $V \in \mathcal{V}_a$ (note that in this case $\overline{\lambda}=a$ does not depend on $\underline{\nu}$) satisfying the equation
\begin{equation}\label{eq:MA2}
h\left(x+\frac{a}{2} \nabla V (x)\right)\det\left(I_d+\frac{a}{2}\nabla^{2}_{x}V(x)\right)=e^{-V(x)}h(x),
\end{equation}
for $\mu$-almost every $x\in \mathrm{dom}(V)=\{x \in \R^d : V(x)<\infty\}.$
Observe that $V=0$ (corresponding to $\underline{\nu}=\mu$) is always solution of \eqref{eq:MA2}.
Let us set
\[
A(\mu)=\inf\{ a\geq 0: V=0 \mbox{  is the only } V \in \mathcal{V}_a \mbox{ s.t.}\int e^{-V(x)}\,d\mu =1 \mbox{  and } \eqref{eq:MA2}  \mbox{ holds } \mu \mbox{ a.s.}\}.
\]
Clearly,
\[
C_{T_2}(\mu) \leq A(\mu).
\]
Moreover, the proof of Proposition \ref{prop:WI} actually shows that
\[
A(\mu) \leq 2\sqrt{C_{W_2I}(\mu)}.
\]
Finally, according to \cite{OV00}, if in addition $\mu$ is assumed to be log-concave, then
\[
\frac{1}{4}C_{LSI}(\mu) \leq C_{T_2}(\mu).
\]
Summarizing the discussion, we get the following result.
\begin{prop}Let $\mu$ be a compactly supported probability measure on $\R^d$ absolutely continuous with respect to Lebesgue measure. Then,
\[
C_{T_2}(\mu)\leq A(\mu)\leq 2\sqrt{C_{W_2I}(\mu)} \leq 4 C_{LSI}(\mu).
\]
If in addition, $\mu$ is log-concave, then
\[
\frac{1}{4}C_{LSI}(\mu) \leq C_{T_2}(\mu)\leq A(\mu)\leq 2\sqrt{C_{W_2I}(\mu)} \leq 4 C_{LSI}(\mu).
\]
\end{prop}
Note that the assumption that the support of $\mu$ is compact can be removed (using Theorem \ref{thm:attainment3} instead of Proposition \ref{prop:attainment} in the proof of Proposition \ref{prop:WI} and in the discussion above). Details are left to the reader.

\smallskip

The interesting conclusion of the result above is that at least in the log-concave case the constant $A(\mu)$ is equal up to universal factors to the best constants in the Talagrand and in the logarithmic Sobolev inequalities. These constants can thus be interpreted in terms of a uniqueness property of a certain Monge-Amp\`ere equation. An interesting question, that will perhaps be discussed elsewhere, would be to try to estimate directly $A(\mu)$ in this framework.

\section{Attainment of the minimum under a weak integrability condition} \label{Sec:attainment}
In this section, we study the existence of a minimizer of the specific functional $F_a$ given by
\begin{equation}\label{eq:EntW2-bis}
F_a(\nu) = aH(\nu|\mu) - \mathcal{T}_c(\nu,\mu),\qquad \forall \nu \in \mathcal{P}_\mu(\X),
\end{equation}
(corresponding to $\alpha(x)=\beta(x)=x$) where $c:\X^2 \to [0,\infty)$ is a power type cost function (as defined in Lemma \ref{lem:power-type}) satisfying the concentration of measure inequality \eqref{eq:concentration} for some $a'>0$ and $r_o\geq0$ that we restate for the reader's convenience: for all $A \subset \X$ such that $\mu(A) \geq 1/2$, it holds
\begin{equation}\label{eq:concentration-bis}
\mu(A_r) \geq 1-e^{-(r-r_o)^{p_o}/a'},\qquad \forall r\geq r_o,
\end{equation}
where $p_o$ is the exponent of $c$ defined in Lemma \ref{lem:power-type}, $A_r = \{ x \in \X: \exists y \in A, \tilde{d}(x,y) \leq r \}$ and $\tilde{d}(x,y) = c^{1/p_o}(x,y)$, $x,y \in \X$.

Our goal is to prove Theorem \ref{thm:attainment3}.
Our strategy will be  based on truncating the cost function, by considering for all positive integer $n$,
\[
\mathcal{T}_{c,n}(\nu_1,\nu_2)= \inf_{\pi} \iint c_n(x,y)\,\pi(dxdy),\qquad \nu_1,\nu_2 \in \mathcal{P}(\X),
\]
where the infimum runs over all couplings between $\nu_1$ and $\nu_2$ and
\[
c_n(x,y)=c(x,y)\wedge n, \qquad x,y\in \X.
\]
Then, introducing (for $a>a'$)
\[
F_{a,n}(\nu) = aH(\nu|\mu) - \mathcal{T}_{c,n}(\nu,\mu), \qquad \nu \in \mathcal{P}_\mu(\X),
\]
we obtain from Proposition \ref{prop:attainment}  that $F_{a,n}$ reaches its minimum at some point $\nu_n$. The rest of the proof consists in showing that $\nu_n$ admits a subsequence
converging to some minimizer of $F_a.$
Let us begin with a simple lemma gathering some properties of the $\nu_n$.

\begin{lem}\label{lem:prop-nu_n} Under the assumptions of Theorem \ref{thm:attainment3}, let $a>a'$.
\begin{enumerate}
\item The sequence $(F_{a,n}(\nu_n))_{n\geq 1}$ converges to $m:=\inf_{\nu \in \mathcal{P}_\mu(\R^d)} F_a(\nu)$. Moreover, the sequence $(\nu_n)_{n\geq1}$ is precompact for the weak topology.
\item For all $n\geq 1$, there exists a couple of continuous functions $(\psi_n,\varphi_n)$ which are $c_n$-conjugate in the sense that $\psi_n(x) = \inf_{y\in \X}\{ -\varphi_n(y) + c_n(x,y)\}$, $x \in \X$ and $\varphi_n(x) = \inf_{x\in \X}\{ -\psi_n(x) + c_n(x,y)\}$, $y \in \X$, such that
\[
\frac{d\nu_n}{d\mu}(x) = \exp\left( \frac{1}{a}\psi_n(x)\right)
\]
and
\[
F_{a,n}(\nu_n) = -\int \varphi_n\,d\mu - a\log\left(\int e^{\frac{1}{a}\psi_n}\,d\mu\right).
\]
\item If there exists two bounded sequences of points $(x_n)_{n\geq 1}$ and $(y_n)_{n\geq 1}$ in $\X$ such that for all $n\geq 1$
\begin{equation}\label{eq:subdiff}
\varphi_n(z) \leq \varphi_n(y_n) + c_n(x_n,z)-c_n(x_n,y_n),\,\forall z\in\X,
\end{equation}
then $(\nu_n)_{n\geq 1}$ admits a subsequence converging to a minimizer of $F_a$.
\end{enumerate}
\end{lem}
\proof\ \\
(1) First let us show that
\begin{equation}\label{eq:limW_2}
\mathcal{T}_{c,n}(\nu,\mu) \to \mathcal{T}_{c}(\nu,\mu),\qquad \text{as } n \to \infty.
\end{equation}
By definition $\mathcal{T}_{c,n}(\nu,\mu) \leq \mathcal{T}_{c}(\nu,\mu)$, which shows that $\limsup_{n \to \infty} \mathcal{T}_{c,n}(\nu,\mu) \leq \mathcal{T}_{c}(\nu,\mu).$ On the other hand, generalizing the argument yielding to \eqref{eq:liminf-utile} in the proof of Theorem \ref{thm:OV}, one gets $\liminf_{n \to \infty} \mathcal{T}_{c,n}(\nu,\mu) \geq \mathcal{T}_{c}(\nu,\mu)$, which gives the desired convergence  \eqref{eq:limW_2}. Now, for all $\nu \in \mathcal{P}_\mu(\X)$,
\[
m \leq F_a(\nu_n)\leq F_{a,n}(\nu_n) \leq F_{a,n}(\nu).
\]
Therefore, letting $n \to \infty$ and using \eqref{eq:limW_2}, one gets
\[
m \leq \liminf_{n\to \infty} F_{a,n}(\nu_n) \leq \limsup_{n\to \infty} F_{a,n}(\nu_n) \leq \lim_{n\to \infty} F_{a,n}(\nu) = F_a(\nu)
\]
Optimizing over $\nu \in \mathcal{P}_\mu(\X)$ completes the proof of the first claim.
Now since $F_{a} \leq F_{a,n}$ and $F_{a,n}(\nu_n)$ is bounded one sees that there exists some $r>0$, such that $\nu_n \in \{F_a \leq r\}$ for all $n\geq 1.$  According to Proposition \ref{prop:bounded-below-refined}, the set $\{F_a \leq r\}$ is  therefore precompact for the weak topology. This implies that the sequence $(\nu_n)_{n\geq 1}$ is itself precompact.

\noindent (2) According to Theorem \ref{thm:CharactMinimizer}, there exists a Kantorovich potential $\psi_n$ for the transport of $\nu_n$ on $\mu$ (for the cost $\mathcal{T}_{c,n}$) such that the density of $\nu_n$ satisfies
\[
\frac{d\nu_n}{d\mu} = \frac{1}{Z_n}\exp\left( \frac{1}{a}\psi_n\right),
\]
with $Z_n = \int e^{\frac{1}{a}\psi_n}\,d\mu$. Letting $\varphi_n(y)=\inf_{x \in \X}\{-\psi_n(x) + c_n(x,y)\}$, it holds $\psi_n(x) = \inf_{y \in \X}\{-\varphi_n(y) + c_n(x,y)\}$.
By definition of Kantorovich potentials, it further holds that
\begin{align*}
\mathcal{T}_{c,n}(\nu_n,\mu) &= \int \psi_n\,d\nu_n + \int \varphi_n\,d\mu.
\end{align*}
On the other hand,
\begin{align*}
aH(\nu_n|\mu) &= \int \psi_n\,d\nu_n - a\log\left(\int e^{\frac{1}{a} \psi_n}\,d\mu\right).
\end{align*}
Thus,
\begin{equation}\label{eq:F(nu_n)}
F_{a,n}(\nu_n) = -\int \varphi_n\,d\mu - a\log\left(\int e^{\frac{1}{a}\psi_n}\,d\mu\right).
\end{equation}
\noindent (3) Let $\varepsilon>0$ and $x_o \in \X$ be an arbitrary point.  Using \eqref{eq:F(nu_n)} at the second line, the inequality $\psi_n(x) \leq -\varphi_n(y_n)+c(x,y_n)$ at the third and \eqref{eq:subdiff} at the last line, one gets
\begin{align*}
&\int e^{\varepsilon c(x,x_o)}\, \nu_n(dx) = \int e^{\varepsilon c(x,x_o) + \frac{1}{a}\psi_n(x)}\,\mu(dx) e^{-\log \left(\int e^{\frac{\psi_n}{a}}\,d\mu\right)}\\
& =  \int e^{\varepsilon c(x,x_o) + \frac{1}{a}\psi_n(x)}\,\mu(dx) \exp\left( \frac{1}{a}F_{a,n}(\nu_n) +\frac{1}{a}\int\varphi_n\,d\mu\right)\\
& \leq \int e^{\varepsilon c(x,x_o) +\frac{1}{a}c_n(x,y_n)}\,\mu(dx) \exp\left( \frac{1}{a}F_{a,n}(\nu_n) + \frac{1}{a} \int \varphi_n(z)-\varphi_n(y_n)\,\mu(dz)\right)\\
& \leq \int e^{\varepsilon c(x,x_o) +\frac{1}{a}c(x,y_n)}\,\mu(dx) \exp\left( \frac{1}{a}F_{a,n}(\nu_n) + \frac{1}{a} \int c_n(x_n,z)-c_n(x_n,y_n)\,\mu(dz)\right).
\end{align*}
Let us check that the last expression above is bounded uniformly in $n$, if $\varepsilon$ is chosen small enough. Indeed, according to Lemma \ref{lem:power-type}, $c=\tilde{d}^{p_o}$ for some metric $\tilde{d}$ and $p_o\geq 1$. Thus, using the triangle inequality for $\tilde{d}$ and the convexity of $t \mapsto t^{p_o}$, we have for all $t\in (0,1)$
\begin{align*}
\cost(x,y_n)= \tilde{d}^{p_o}(x,y_n)\leq \left(\tilde{d}(x,x_o)+\tilde{d}(x_o,y_n)\right)^{p_o} &= \left( (1-t) \frac{\tilde{d}(x,x_o)}{1-t}+t \frac{\tilde{d}(x_o,y_n)}{t}\right)^{p_o}\\
& \leq \frac{\tilde{d}^{p_o}(x,x_o)}{(1-t)^{p_o-1}}+\frac{\tilde{d}^{p_o}(x_o,y_n)}{t^{p_o-1}}
\end{align*}
According to Item (2) of Lemma \ref{lem:comparison}, the concentration inequality \eqref{eq:concentration-bis} implies that
\[
\int e^{\delta c(x,x_o)}\,\mu(dx)<\infty,
\]
for all $\delta<1/a'.$
Thus, if $\varepsilon$ and $t$ are chosen so that $\varepsilon + \frac{1}{a(1-t)^{p_o-1}} <1/a'$, we have
\[
\sup_{n\geq 1}\int e^{\varepsilon c(x,x_o) + \frac{1}{a}c(x,y_n)}\,\mu(dx) <+\infty.
\]
Similarly, $\int c_n(x_n,z)\,\mu(dz) \leq \int c(x_n,z)\,\mu(dz) \leq 2^{p_o-1} \int c(x_o,z)\,\mu(dz) + 2^{p_o-1} c(x_n,x_o)$. Since $x_n,y_n$ are bounded, $F_{a,n}(\nu_n)$ converges and $c_n\geq0$, we conclude that
\[
\sup_{n\geq 1}\exp\left( \frac{1}{a}F_{a,n}(\nu_n) + \frac{1}{a} \int c_n(x_n,z)-c_n(x_n,y_n)\,\mu(dz)\right)<+\infty.
\]
In conclusion, if $\varepsilon$ is small enough, one has
\begin{equation}\label{eq:UInu_n}
\sup_{n\geq 1} \int e^{\varepsilon c(x,x_o)}\, \nu_n(dx) <+\infty.
\end{equation}
According to Item (1), the sequence $\nu_n$ is precompact, and so it admits a subsequence (also denoted by $\nu_{n}$ for simplicity) converging weakly to some $\underline{\nu} \in \mathcal{P}(\X).$ Moreover, it is easy to see from \eqref{eq:UInu_n} that
\[
\lim_{k\to \infty} \sup_{n\geq 1}\int c(x,x_o)\mathbf{1}_{c(x,x_o)\geq k}\,\nu_n(dx) = 0.
\]
By \cite[Theorem 6.9]{Vil09}, the convergence
\[
\widetilde{W}_{p_o}(\nu_{n},\underline{\nu}) \to 0,\qquad \text{as } n\to\infty.
\]
 also holds  true for  $\widetilde{W}_{p_o}$ the Wasserstein distance associated with the metric $\tilde{d}$.  Therefore, it also holds
\[
\mathcal{T}_c(\nu_n,\mu) = W_{p_o}^{p_o}(\nu_n,\mu) \to W_{p_o}^{p_o}(\underline{\nu},\mu) = \mathcal{T}_c(\underline{\nu},\mu),
\]
as $n \to \infty.$
Together with the lower semicontinuity of $H(\,\cdot\,|\mu)$, this  immediately implies that
\[
F_a(\underline{\nu}) \leq \liminf_{n\to \infty}F_a(\nu_n)
\]
Since $F_{a}(\nu_n) \leq F_{a,n}(\nu_n)$ and, according to Item (1), $F_{a,n}(\nu_n) \to m=\inf_{\nu \in \mathcal{P}_\mu(\X)} F_a(\nu)$, as $n \to \infty$, we conclude that $F_a(\underline{\nu}) \leq m$, and so $F_{a}(\underline{\nu})=m$, which completes the proof.
\endproof
\begin{rem}
The reader familiar with the notion of $\Gamma$-convergence will have noticed that some steps in the preceding proof could be derived from general principles available for instance in the classical text book \cite{DM93}.
Let us emphasize some simplifications that can be performed using tools from \cite{DM93}.
First of all, as shown in Point (1), the sequence of functions $F_{a,n}$ converges pointwise to $F_a$. Since this sequence is non-increasing, it follows from \cite[Proposition 5.7]{DM93} that $F_{a,n}$ converges to $\mathrm{sc}^- (F_a)$ in the sense of $\Gamma$-convergence. Here, by $\mathrm{sc}^- (F_a)$, we denote the lower semicontinuous envelop of $F_a$, that is to say the greatest lower semicontinuous function below $F_a$ (see \textit{e.g.} \cite[Chapter 3]{DM93}). On the other hand, as shown again in Point (1), the sequence $\nu_n$ is precompact. Therefore, it follows from \cite[Theorem 7.4]{DM93}, that the function $\mathrm{sc}^-(F_a)$ attains its minimum (which also follows from the fact that, as shows Proposition \ref{prop:bounded-below-refined}, $F_a$ has precompact sublevel sets - is \emph{coercive} in the terminology of \cite{DM93}- and from \cite[Point (b) of Theorem 3.8]{DM93}) and
\[
F_{a,n}(\nu_n) \to \min \mathrm{sc}^- (F_a).
\]
Now, according to \cite[Point (c) of Theorem 3.8]{DM93}, $\min \mathrm{sc}^- (F_a) = \inf F_a.$ According to \cite[Corollary 7.20]{DM93}, it also follows that if $\underline{\nu}$ is any cluster point of $\nu_n$, then $\underline{\nu}$ is a minimizer of $\mathrm{sc}^-(F_a).$ Therefore, it holds $\inf F_a=\mathrm{sc}^-F_a(\underline{\nu})\leq F_a(\underline{\nu}).$
So if one can show that  $\nu_n$ admits a cluster point $\underline{\nu}$ which satisfies $\mathrm{sc}^-(F_a)(\underline{\nu}) = F_a(\underline{\nu})$, then $\underline{\nu}$ will be a minimizer of $F_a.$ This is what we prove in Point (3) of Lemma \ref{lem:prop-nu_n} and in the rest of the proof of Theorem \ref{thm:attainment3} below.
\end{rem}

Now, the question is to prove the existence of bounded sequences $x_n,y_n$ as in Item (3) of Lemma \ref{lem:prop-nu_n}.
We begin by  stating a classical lemma showing that the weak convergence of a sequence of probability measures implies the convergence of their supports in the sense of Kuratowski.
\begin{lem}\label{lem:support}
Let $(\gamma_n)_{n\geq 1}$ be a sequence of probability measures defined on some Polish space $(E,d)$ converging weakly to some probability measure $\gamma.$ Then for all point $z$ in the support of $\gamma$, there exists a sequence $(z_n)_{n\geq 1}$ such that for all $n\geq 1$ $z_n$ belongs to the support of $\gamma_n$ and $z_n$ converges to $z$ as $n$ tends to $\infty.$
\end{lem}
This result is proved in \textit{e.g.} \cite[Proposition 5.1.8]{AmbGigSav-08}.

Now we are ready to complete the proof of Theorem \ref{thm:attainment3}.
\proof[Proof of Theorem \ref{thm:attainment3}]
Here we use the notations introduced in Lemma \ref{lem:prop-nu_n}.
Let $(\nu_n)_{n \geq 1}$ be a sequence of minimizers of $F_{a,n}$, for some $a>a'.$ According to Item (1) of Lemma \ref{lem:prop-nu_n}, the sequence $(\nu_n)_{n\geq 1}$ is precompact for the weak topology.
Therefore, one can assume without loss of generality that $\nu_n$ converges to some $\nu \in \mathcal{P}_\mu(\X).$ For all $n\geq 1$, let $\pi_n$ be an optimal coupling between $\nu_n$ and $\mu$ for the cost $\mathcal{T}_{c,n}.$ Since the marginals of $\pi_n$ are converging, a classical argument shows that $\pi_n$ is a tight sequence (see \textit{e.g.} \cite[Theorem 4.4]{Vil09}), and so according to Prokhorov Theorem, it admits at least one converging subsequence, still denoted by $\pi_n$ in the sequel. Take an arbitrary point $\bar{z}=(\bar{x},\bar{y})$ in the support of $\pi$ ;  according to Lemma \ref{lem:support} above, there exists a sequence of points $z_n=(x_n,y_n)$ such that $z_n$ belongs to the support of $\pi_n$ and  $z_n \to z$  as $n \to \infty.$
For all $n\geq 1$, since $\pi_n$ and $\psi_n,\varphi_n$ are optimal, it holds
\[
\int c_n(x,y)-(\psi_n(x)+\varphi_n(y))\,\pi_n(dxdy) =0.
\]
By definition, $\psi_n(x)=\inf_{y\in \X}\{-\varphi_n(y) + \tilde{d}(x,y)^{p_o} \wedge n\}$, $x \in \X$, where $\tilde{d}$ and $p_o$ have been introduced in Lemma \ref{lem:power-type}. It is not difficult to check that for any $y \in \X$, the function $x\mapsto \tilde{d}(x,y)^{p_o} \wedge n$ is $p_on^{(p_o-1)/p_o}$-Lipschitz with respect to $\tilde{d}$. As an infimum of Lipschitz functions, $\psi_n$ is also $p_on^{(p_o-1)/p_o}$-Lipschitz with respect to $\tilde{d}$, from which one deduces easily that $\psi_n$ is continuous on $\X.$ The same argument applies to $\varphi_n$, so the integrand in the integral above is continuous and non-negative, and thus $\psi_n(x)+\varphi_n(y) = c_n(x,y)$ for all $(x,y)$ belonging to the support of $\pi_n$.
In particular, $\psi_n(x_n)+\varphi_n(y_n)=c_n(x_n,y_n)$ and since $\psi_n(x_n) = \inf_{z \in \X }\{-\varphi_n(z) + c_n(x,z)\}$ one concludes that
\[
\varphi_n(z) \leq \varphi_n(y_n) + c_n(x_n,z) - c_n(x_n,y_n),
\]
for all $z \in \X.$ Using Item (3) of Lemma \ref{lem:prop-nu_n}, one concludes that $F_a$ admits a minimizer.
\endproof

\section{Links with the characterization of moment measures}\label{Sec:Moment-Measures}
As mentioned in the introduction, Equation \eqref{eq:twistedMM} and the minimization problem of the functional $F_a$ given in \eqref{eq:EntW2-bis} in the case $c=d^2$ feature close connections with the recent work \cite{CEK15} by Cordero-Erausquin and Klartag on the characterization of moment measures. Let us recall that a Borel measure $\mu$ on $\er^d$ is said to be a moment measure for a (convex) function $\phi:\R^d \to \R\cup\{+\infty\}$ such that $0<\int e^{-\phi}\,dx<+\infty$ if $\nabla \phi$ pushes forward the measure $\nu_\phi(dx)=\frac{e^{-\phi(x)}}{\int e^{-\phi(x)}\,dx}\,dx$ towards $\mu$; or equivalently $\nabla \phi_\#\nu_\phi=\mu$. This notion of moment measures finds applications in differential geometry, partial differential equations (in particular Monge-Amp\`ere equation) and the study of log-concave measures (we refer the interested reader to \cite{CEK15} and references therein for more details).

In \cite{CEK15}, the authors obtain a new characterization of moment measures, showing that \emph{any} Borel measure $\mu$ with positive finite total mass on $\er^d$, such that its support has dimension $d$ and with $0$ as barycenter can be represented as a moment measure for some (unique up to translations) essentially continuous convex function $\phi.$ This characterization was obtained though the study of the well posedness of maximizers of the functional
\[
I_{\mu}(\phi)=\log\left(\int e^{-\phi^*}\,dx\right)-\int \phi\,d\mu
\]
on the space of (proper) convex functions $\phi$ such that $0<\int e^{-\phi^*(x)}\,dx<+\infty$, where
\[
\phi^*(x)=\sup_{x\in \R^d}\{x\cdot y - \phi(y)\}, \qquad x \in \R^d
\]
is the Legendre-Fenchel transform of $\phi$. Under the assumptions on $\mu$ recalled above, the authors show that the maximum value of $I_{\mu}$ is attained on a $\mu$-integrable convex function $\overline{\phi}$ (unique up to translations) such that $\mu= \nabla \overline{\phi}_\#\nu_{\overline{\phi}}$. The  arguments for  existence, uniqueness and characterization of the maximizer of $I_{\mu}$ rely mostly on convex analysis and functional inequalities.

More recently, Santambrogio \cite{S15} provides a dual counter-part of the results obtained in \cite{CEK15}, considering the minimization, over the space $\Pp_1(\er^d)$ (of probability measures admitting a finite first moment), of
\[
J_{\mu}(\nu)=H(\nu\,|\,\mbox{Leb})+\sup_{\pi\in\Pi(\mu,\nu)}\int x\cdot y\,\pi(dxdy),
\]
where $H(\nu\,|\,\mbox{Leb})$ is minus the Shannon entropy of $\nu$:
\begin{equation*}
H(\nu\,|\,\mbox{Leb})=\left\{
\begin{aligned}
&\int \log\left(\frac{d\nu}{dx}\right)\,\nu(dx)\,\mbox{whenever the density }\frac{d\nu}{dx}\,\mbox{exists},\\
&+\infty\,\mbox{otherwise}.
\end{aligned}
\right.
\end{equation*}
The minimization problem related to $I_{\mu}$ is dual to the optimization problem in \cite{CEK15} in the sense that, as shown in Section 6 of \cite{S15},
\begin{equation}
\sup_{\phi} I_{\mu}(\phi)=-\inf_{\nu\in\Pp_1(\R^d)} J_{\mu}(\nu)\label{eq:DualMM}.
\end{equation}
Assuming that $\mu$ still satisfies the hypotheses of \cite{CEK15}, existence and uniqueness of a minimizer $\underline{\nu}$ in $\Pp_1(\R^d)$ of $J_\mu$ rely on optimal transport theory and coupling techniques and, using the sub-differential calculus on $\Pp(\er^d)$ (in a similar way as  we did in Theorem \ref{thm:CharactMinimizer}), $\underline{\nu}$ is characterized by the property that
\[
\log\left(\frac{d\underline{\nu}}{dx}\right)=-\phi,
\]
where $\phi$ is a convex function whose gradient is the $\Tt_2$-optimal map pushing forward $\underline{\nu}$ to $\mu$. This provides an alternative (and in some sense more direct) characterization of $\mu$ as a moment measure. Note that the sub-differential calculus mentioned above is understood in its usual sense, and is not in the Wasserstein sense of Ambrosio-Gigli-Savar\'e \cite{AmbGigSav-08, AGS14} (see also \cite{Santambook}).

The minimization problem of $F_a$ given by \eqref{eq:EntW2-bis} with $c(x,y)=\frac{1}{2}|x-y|^2$ and  $a=1$ can be related to the minimization problem of $J_\mu$ since, if $\mu\in\Pp_2$,
\begin{equation}\label{eq:LinkS15}
F_a(\nu)=J_{\mu}(\nu)-\int \log\left(\frac{d\mu}{d\gamma}\right)(x)\,\nu(dx)-\int \frac{|x|^2}{2}\,\mu(dx),
\end{equation}
where $\gamma$ denotes a centered Gaussian distribution on $\R^d$ with covariance matrix $\mathrm{Id}$. In particular, both problems coincide in the case $\mu=\gamma$.

In the same spirit as \eqref{eq:DualMM}, we can give a dual formulation of our minimization problems:
\begin{prop}
Let $\mu$ be a probability measure on $\X$ and consider the functional $F_a(\nu) = aH(\nu|\mu) - \mathcal{T}_c(\nu,\mu)$, $\nu\in \mathcal{P}_\mu(\X)$, where $c:\X^2\to\R^+$ is some continuous cost function such that $\iint e^{\delta c(x,y)}\,\mu(dx)\mu(dy) <\infty$ for some $\delta\geq 1/a.$
Then,
\[
\inf_{\nu \in \mathcal{P}_\mu(\X)} F_a(\nu) =  \inf_{(\psi,\varphi) \in \mathcal{F}} \left\{ - \int\varphi\,d\mu -a \log \int e^{\psi/a}\,d\mu \right\},
\]
where $\mathcal{F}$ denotes the set of couples of bounded continuous functions $(\psi,\varphi)$ such that $\psi(x) + \varphi(y) \leq c(x,y)$, for all $x,y \in \X.$
\end{prop}
This result is essentially a rewriting of Bobkov-G\"otze dual formulation of transport-entropy inequalities \cite{BG99} (see \cite[Section 3]{SurveyGL} for general statements).
\proof
According to Kantorovich's duality, it holds
\[
\mathcal{T}_c(\nu,\mu)  = \sup_{(\psi,\varphi)\in \mathcal{F}} \left\{ \int \psi(x)\,\nu(dx) + \int \varphi(y)\,\mu(dy)\right\},
\]
so that
\begin{align*}
\inf_{\nu \in \mathcal{P}_\mu(\X)} F_a(\nu) & = \inf_{\nu \in \mathcal{P}_\mu(\X)}\inf_{(\psi,\varphi) \in \mathcal{F}} \left\{ aH(\nu|\mu) - \int \psi\,d\nu -\int \varphi\,d\mu\right\}\\
& = \inf_{(\psi,\varphi) \in \mathcal{F}} \left\{ - \int\varphi\,d\mu + \inf_{\nu \in \mathcal{P}_\mu(\X)}\left\{ aH(\nu|\mu) - \int \psi\,d\mu\right\} \right\}.
\end{align*}
According to a well known duality formula for the relative entropy,
\[
\sup_{\nu \in \mathcal{P}_\mu(\X)} \left\{\int f\,d\nu - H(\nu|\mu)\right\} = \log \int e^{f}\,d\mu.
\]
Therefore,
\[
\inf_{\nu \in \mathcal{P}_\mu(\X)} F_a(\nu) = \inf_{(\psi,\varphi) \in \mathcal{F}} \left\{ - \int\varphi\,d\mu -a \log \int e^{\psi/a}\,d\mu \right\}
\]
\endproof
In particular, if $c(x,y)=\frac{1}{2}|x-y|^2$, $x,y \in \R^d$, then $(\psi,\varphi) \in \mathcal{F}$ if and only if $f(x) = -\psi(x)+|x|^2/2$ and $g(y)=-\varphi(y)+|y|^2/2$ satisfy $f(x)+g(y) \geq x\cdot y$, for all $x,y \in \R^d.$ From this it is not difficult to see that
\begin{align*}
\inf_{\nu\in\mathcal{P}_\mu(\X)} F(\nu)=-\sup_{\phi}\left\{a\log\left(\int e^{-\phi/a}e^{\frac{|x|^2}{2a}}d\mu\right)-\int \phi^*\,\mu(dx)\right\}-\int \frac{|x|^2}{2}\,\mu(dx),
\end{align*}
where the supremum applies over all convex $L^1(\mu)$-function and where $\phi^*$ is the Legendre-Fenchel transform of $\phi$. Details are left to the reader.

Hence, at first sight, the minimization of the functional \eqref{eq:EntW2-bis} might be considered using the general techniques used in \cite{S15} or \cite{CEK15}.
Yet, to obtain existence of a minimizer, a direct adaption of the proofs in those works would    require stronger assumptions on $\mu$ than ours, in order to ensure that the minimum of $F_a$ is attained (for instance, the  uniform integrability of the second moments of the sequence of minimizers of $F_a$  used in \cite{S15} would require $\mu$  to have some exponential moments of order strictly larger that $2$). Therefore, the truncation technique used to construct a minimizer to \eqref{eq:EntW2-bis} provides an approach alternative to  \cite{CEK15} and \cite{S15}, which furthermore can be extended to more general settings. This will be the subject of future works.

Let us close this section with a remark on the problem of  uniqueness of a minimizer of \eqref{eq:EntW2-bis} and another distinction between our problem and the problem in \cite{S15} and \cite{CEK15}. Uniqueness of a minimizer to $\inf J_{\mu}$ follows directly from the (strict) displacement convexity of $J_{\mu}$ with respect to $W_2$-geodesics (see \textit{e.g.} \cite{AmbGigSav-08} and \cite{Vil09}), while the uniqueness of a maximizer to $\sup I_{\mu}$ is obtained from Pr\'ekopa's inequality (these two 'convexity' properties being dual of each other).
In our setting, using similar ingredients, one can easily prove that the functional $F_a$ associated to a quadratic cost is strictly geodesically convex, when the reference probability measure $\mu$ is uniformly log-concave. More precisely one has the following.
\begin{prop}
Suppose that $\mu(dx)=e^{-V(x)}\,dx$ with $V:\R^d \to \R$ a function of class $\mathcal{C}^2$ such that $\mathrm{Hess}\,V \geq K\mathrm{Id}$, for some $K>0.$ Then for all $a> 2/K$, the functional $F_a$ of \eqref{eq:EntW2-bis} with $c(x,y)=|x-y|^2$, $x,y\in \R^d$ is strictly geodesically convex.
\end{prop}
\proof
According to \cite[Theorem 7.3.2]{AmbGigSav-08}, if $\nu_0,\nu_1 \in \mathcal{P}_2(\R^d)$ (the space of probability measures having finite second moments), it holds, for all $t \in [0,1]$
\[
W_2^2(\nu_t,\mu) \geq (1-t)W_2^2(\nu_0,\mu)+t W_2^2(\nu_1,\mu)-t(1-t)W_2^2(\nu_0,\nu_1),
\]
for all constant speed geodesic $(\nu_t)_{t\in [0,1]}$ (for the $W_2$ metric) joining $\nu_0$ to $\nu_1$.
On the other hand, according to \textit{e.g.} \cite[Theorem 17.15]{Vil09}, the relative entropy functional satisfies, for all $\nu_0,\nu_1 \in \mathcal{P}_\mu(\R^d),$ and all constant speed geodesic $\nu_t$ joining $\nu_0$ to $\nu_1$,
\begin{equation}\label{eq:disp-conv}
H(\nu_t |\mu) \leq (1-t)H(\nu_0|\mu) + t H(\nu_1|\mu) - \frac{K}{2} t(1-t) W_2^2(\nu_0,\nu_1).
\end{equation}
So it follows immediately, that for $a >2/K$, $F_a$ satisfies
\[
F_a(\nu_t) < (1-t)F_a(\nu_0)+tF_a(\nu_1),\qquad \forall \nu_0\neq \nu_1  \in \mathcal{P}_\mu(\R^d),\qquad \forall t\in (0,1),
\]
which completes the proof.
\endproof
Of course, strict convexity ensures uniqueness of the minimizer. But the assumption of uniform convexity of the potential $V$ is too strong to be really interesting for our purpose. Indeed under this assumption, \eqref{eq:disp-conv} immediately implies that $\mu$ satisfies $\T_2(2/K)$. Let us recall this well known argument. Taking $\nu_0=\mu$ and using that $H(\nu_t|\mu)\geq 0$, one immediately gets from \eqref{eq:disp-conv} that $tH(\nu_1|\mu) - \frac{K}{2}t(1-t)W_2^2(\nu_1,\mu)\geq0$. Dividing by $t$ and then letting $t\to1$ proves the claim. Nevertheless, studying the uniqueness of a minimizer to $F_a$ is an interesting question which needs to be handled by suitable techniques and which will be deepened in future works.
\bibliographystyle{plain}
\bibliography{Bib-FGJ}
\end{document}